\documentclass{amsart}
\usepackage{amssymb,stmaryrd}
\usepackage{amsfonts}
\usepackage{amstext}
\usepackage{algorithmic}
\usepackage{algorithm}
\usepackage{graphicx}
\usepackage{epstopdf}
\usepackage[all]{xy}
\usepackage{MnSymbol}
\usepackage{tikz}

\numberwithin{equation}{section}
\newtheorem{theorem}{Theorem}[section]
\newtheorem{lemma}[theorem]{Lemma}
\newtheorem{proposition}[theorem]{Proposition}

\theoremstyle{definition}

\theoremstyle{remark}

\newcommand{\R}{{\mathbb{R}}}

\newcommand{\C}{{\mathbb{C}}}
\newcommand{\Z}{{\mathbb{Z}}}

\newcommand{\F}{{\mathbb{F}}}

\newcommand{\<}{{\langle}}

\renewcommand{\)}{{)}}
\renewcommand{\>}{{\rangle}}

\newcommand{\CN}{{\mathcal{N}}}
\newcommand{\wedgeq}{{\wedge\kern-5pt\cdot}}

\renewcommand{\ker}{{\rm{ker}}}

\newcommand{\tens}{\otimes}

\newcommand{\id}{{\rm id}}

\newcommand{\extd}{{\rm d}}
\newcommand{\del}{{\partial}}
\newcommand{\eps}{\epsilon}

\newcommand{\twoForm}[2]{e^{#1} \otimes e^{#2}}
\newcommand{\sigFun}[2]{ \sigma(\twoForm{#1}{#2}) }

\begin{document}

\title{Quantum geometry of Boolean algebras and de Morgan duality}


\author{S. Majid}
\address{Queen Mary University of London\\
School of Mathematics, Mile End Rd, London E1 4NS, UK}

\email{s.majid@qmul.ac.uk}
\date{}

\begin{abstract} We take a fresh look at the geometrization of logic using the recently developed tools of `quantum Riemannian geometry' applied in the digital case over the field $\F_2=\{0,1\}$, extending de Morgan duality to this context of differential forms and connections. The 1-forms correspond to graphs and the exterior derivative of a subset amounts to the arrows that cross between the set and its complement. The line graph $0-1-2$ has a non-flat but Ricci flat quantum Riemannian geometry. The previously known four quantum geometries on the triangle graph, of which one is curved, are revisited in terms of left-invariant differentials, as are the quantum geometries on the dual Hopf algebra, the  group algebra of $\Z_3$. For the square,  we find a moduli of four quantum Riemannian geometries, all flat, while for an $n$-gon with $n>4$ we find a unique one, again flat. We also propose an extension of de Morgan duality to general algebras and differentials over $\F_2$.  \end{abstract}
\maketitle 

\section{Introduction}

This is the fifth in a series of papers in which quantum groups and quantum geometry are developed over the field $\F_2$ of two elements, the `digital case'\cite{BasMa,MaPac1,MaPac, MaPacd}. Particularly,  \cite{MaPac} classified low-dimensional parallelisable quantum Riemannian geometries over $\F_2$ with a unique basic top form in degree 2,  and \cite{MaPacd} classified low-dimensional Hopf algebras and bialgebras, and quasitriangular structures on the former. The particular style of noncommutative geometry here is a `constructive approach'  coming out of (but not limited to) experience with the geometry of quantum groups and quantum spacetime models, including the first model \cite{MaRue} of the latter with quantum group symmetry. This approach is somewhat different from Connes' approach to noncommutative geometry\cite{Con} founded in cyclic cohomology and spectral triples or `Dirac operators' but not incompatible with it \cite{BegMa2}. In recent years, it was developed particularly with bimodule connections\cite{DV2} in a series of works with Beggs, as now covered in the book \cite{BegMa}. See also some of the recent literature such as \cite{BegMa1, Ma:gra, Ma:haw, Ma:squ, MaPac, ArgMa}.  One of the key features of this approach is that it can be explored for any unital algebra $A$ over any field. More details are in the preliminaries Section~\ref{secpre}, but the bare essentials are that we start with a first order calculus $(\Omega^1,\extd)$ as a bimodule of 1-forms and a differential from functions. A metric is then an element in $\Omega^1\tens_A\Omega^1$ with certain properties such as the existence of a bimodule map inverse $(\ ,\ ):\Omega^1\tens_A\Omega^1\to A$, and a linear connection is a map $\nabla:\Omega^1\to\Omega^1\tens_A\Omega^1$ with certain properties. This is called a `quantum Levi-Civita' connection (QLC) when it is torsion free and metric compatible. Also of interest is a more symmetric notion of `weak quantum Levi-Civita' connection (WQLC) when the torsion and a certain cotorsion vanish. 

In the present paper, we apply this formalism of quantum Riemannian geometry\cite{BegMa} particularly in the digital discrete geometry case of the algebra $A=\F_2(X)$ for a set $X$. In this case the choice of differential structure in the sense of 1-forms $\Omega^1$ amounts to a graph on $X$. Developing quantum differential geometry on $X$, over any field, therefore includes the potential for a new generation of geometric invariants of graphs\cite{Ma:gra}, with the $\F_2$ case being our focus now. We will see that the metric in the $\F_2$ case is unique, so an example of such a graph invariant would be the moduli space of QLCs for this metric, as well as the sub-moduli of flat ones, or ones with conditions on the curvature such as Ricci flatness. The moduli of WQLCs will also be interesting in this respect, as well as useful to compute as an intermediate step. This kind of moduli space analysis for connected graphs on small sets $X$ is the topic of Section~4 and one of the main results from the point of view of continuity with previous work. Our results for $|X|=2,3$ are complete for the given choices of $\Omega^2$, while for $|X|>3$ we analyse only the polygon case with natural exterior algebras coming from a Cayley graph point of view. A full analysis for the 6 connected graphs for $|X|=4$ in the same spirit as for $|X|=2,3$ is certainly possible but deferred for further work. The specification of $\Omega^2$ goes potentially beyond the graph data, but one of our goals is to explore the proposal in \cite{BegMa} of successive canonical quotients $\Omega_{med}$ and $\Omega_{min}$ of the maximal prolongation  $\Omega_{max}$ of any graph calculus. More details of the set-up are in Section~\ref{secX}. Our results in Section~\ref{secex} subsume the part of the computer classification in \cite{MaPac} that relates to the algebra of functions on a set with $|X|\le 3$ but without the strong assumptions on the exterior algebra made there, and now using more informative algebraic methods. 

The other main goal of the paper is to view the algebras $\F_2(X)$ as complete atomic Boolean algebras and translate their quantum Riemannian geometry to the power set $P(X)$ of subsets of a set with intersection and union, so as to obtain a geometric picture of de Morgan duality. We find that the formalism reduces in the power set Boolean case to a reasonable theory of Riemannian geometry at the level of Venn diagrams on graphs. The first layer of this is the differential structure, which (as usual for functions on a set) amounts to arrows between the underlying elements together together with, which is new, a noncommutative extension of $\cap$ and $\cup$ to include such arrows. The exterior derivative $\extd a$ of a subset will be the subset of arrows that cross between $a$ and its complement $\bar a$.  We will also exhibit Riemannian connections and curvature in some nontrivial examples. Boolean algebras are also the model for propositional logic where $a\subseteq b$ appears as $a\Rightarrow b$ for entailment and our  results can in principle be translated in these terms also. Translation of the formalism is in Section~\ref{secboo}, followed by Section~\ref{secDM} for the new feature of de Morgan duality. We show that this duality (interchanging a set with its complement and $\cap$ with $\cup$ or $a\Rightarrow b$ to $\bar b\Rightarrow \bar a$ in propositional terms) indeed extends to the quantum Riemannian geometry as, in some sense, an extension of  diffeomorphism invariance. 

What the extension of de Morgan duality to quantum Riemannian geometry and basic features of the latter, such  as connections and curvature, mean for logic is not something we can hope to address here. This would be a direction for further work.  Likewise, although the present work is confined to mathematics and does not discuss physics, there could be a possible philosophical basis whereby de Morgan duality, when suitably extended, may be relevant to quantum gravity \cite{Ma:pri,Ma:qg, Ma:mar}. Suffice it to say that while quantum theory is intuitionistic in character as in a Heyting algebra, where we relax the rule that $a\cup\bar a=$everything,  gravity might be expected to be cointuitionistic in character in the de Morgan dual sense, as in a coHeyting algebra, where we relax the rule that $a\cap\bar a=\emptyset$. The latter has also been proposed for other reasons in \cite{Law}  as geometric in nature with $\del a=a\cap \bar a$ a kind of boundary of $a$.  Thus, a kind of duality between quantum theory and geometry, also linked to Hopf algebra duality and quantum Born reciprocity\cite{Ma:pla}, could speculatively have its primitive origin in something like de Morgan duality.  In this respect our study includes an element of gravity in the loose sense of a typically curved metric and an element of quantum theory in the minimal sense that differential forms on Boolean algebras do not commute with algebra elements. 

$\F_2$-geometry is also interesting in its own right and could have other applications, such as to the transfer of geometric ideas to digital electronics\cite{MaPac1,MaPac,MaPacd} and to quantum computing\cite{Ma:mar}. With this in mind, Section~\ref{secbarA} contains a proposal to extend de Morgan duality to general algebras over $\F_2$. This section also reworks the quantum Riemannian geometry of $\F_2\Z_3$ from \cite{MaPac} in terms of basic invariant 1-forms $e^\pm$ parallel to our new treatment of the triangle case $\F_2(\Z_3)$ in Section~\ref{sectri}, of which it is the Hopf algebra dual model.  The work \cite{MaPac} finds many other commutative algebras with rich quantum Riemannian geometries but which are not Boolean algebras, such as $\F_8$ as a 3-dimensional algebra over $\F_2$, or the 4-dimensional Hopf algebra $A_2$ in \cite{BasMa} as part of a family $A_d$ over any $\F_p$. For interesting noncommutative examples, one could consider  the 4-dimensional noncommutative noncocommutative self-dual Hopf algebra ${\sl dsl}_2$ found in \cite{MaPac1}, albeit its intrinsic  quantum Riemannian geometry remains to be explored. Section~\ref{seccon} provides some concluding remarks. 

In this paper, we will take both an $\F_2$-algebra and a Boolean algebra point of view, with $\bar{\ }$ always denoting complementation or its generalisation (never complex conjugation). 

\section{Preliminaries on quantum Riemannian geometry} \label{secpre}

We will be interested in unital algebras $A$ over $\F_2$, which just amounts to a unital ring such that every element is its own additive inverse. Among these,  Boolean algebras are characterised as having all elements square to themselves. We will be particularly interested in the complete atomic case of power sets where our $\F_2$-algebras are of the from $A=\F_2(X)$ for some set $X$, which could be infinite. In the present section, we very briefly recap the geometric formalism at our disposal, which works over any field. More details can be found in \cite{BegMa}, while our treatment of the discrete set case first appeared in \cite{Ma:gra}. 

\subsection{Generalities} We start with some preliminaries on quantum Riemannian geometry as developed on any unital algebra $A$ over any field $k$, see \cite{BegMa} and references therein such as \cite{DV2,Mou,Ma:rie,Ma:gra,BegMa1}.   

The first step is to define a differential structure on $A$ as an $A$-bimodule $\Omega^1$ of `1-forms' and a map $\extd: A\to \Omega^1$ obeying the bimodule Leibniz rule and such that $A\extd A=\Omega^1$. There is also a reasonable notion of diffeomorphism $\phi:A\to B$ between algebras with differentials namely $\phi$ an algebra map which extends to a uniquely determined $A$-bimodule map $\phi_*$ such that
\begin{equation}\label{diffeo}\begin{array}{rcl} \Omega^1_A & {\buildrel \phi_*\over\longrightarrow} & \Omega^1_B\\ \extd\uparrow & & \uparrow\extd\\
A &{\buildrel\phi\over\longrightarrow}&B\end{array}\end{equation}
commutes, for the respective $\Omega^1$ bimodules as marked with $\Omega^1_B$ an $A$-bimodule by pull-back on $\phi$.  There is a `universal first order calculus' $\Omega^1_{uni}=\ker: A\tens A\to A$ with $\extd_{uni}(a)=1\tens a-a\tens 1$ with the property that any other calculus 
is a quotient (the image of a surjection from $\Omega^1_{uni}$.  A calculus is called {\em connected} if $\ker\extd=k1$. 

We will also need to specify higher differential forms (at least to  degree 2 for our purposes) as a DGA (or differential graded algebra)  $(\Omega,\extd)$ with $A,\Omega^1$ in degree 0,1 respectively and with $\extd$ obeying a graded-Leibniz rule and $\extd^2=0$. We require $\Omega$ to generated by its degree 0,1, in which case one says that $\Omega$ is an exterior algebra. An algebra map $\phi$ is {\em strongly differentiable} if it extends similarly as in degree 1, now to $\phi_*:\Omega_A\to \Omega_B$ as a map of DGA's.  Every $A,\Omega^1$ has a canonical `maximal prolongation' $\Omega_{max}$ with the fewest relations to obtain an exterior algebra, so other choices are quotients of this. The maximal prolongation of $\Omega^1_{uni}$ is a well-known universal exterior algebra $(\Omega_{uni},\extd_{uni})$, where $\Omega_{uni}^n\subset A^{\tens(n+1)}$ is the joint kernel of the adjacent product maps and $\extd$ is given by insertion of 1 summed with signs over the different positions. The product of $\Omega$ is denoted $\wedge$. 

In this context, we define a generalised quantum metric as an element $g\in\Omega^1\tens_A\Omega^1$ which is nondegenerate in the sense of the existence of a bimodule map $(\ ,\ ):\Omega^1_A\tens \Omega^1_A\to A$ such that $((\ ,\ )\tens \id)(\omega\tens g)=\omega=(\id\tens(\ ,\ ))(g\tens\omega)$ for all $\omega$ (i.e., $\Omega^1$ is left and right dual to itself in the monoidal category of $A$-$A$-bimodules). By a quantum metric, we mean this data with a further `quantum symmetry' condition, usually taken to be $\wedge(g)=0$, but there could be more natural choices in specific contexts (such as edge symmetry in the graph case).

The remaining ingredients are connections and curvature. A (left) bimodule connection  on $\Omega^1$ is a map $\nabla:\Omega^1\to \Omega^1\tens_A\Omega^1$ subject to the twin Leibniz rules
\begin{equation}\label{nabderiv} \nabla(a\omega)=\extd a\tens\omega+a\nabla\omega,\quad \nabla(\omega a)=(\nabla \omega)a+\sigma(\omega\tens\extd a),\quad \forall a\in A, \omega\in \Omega^1\end{equation}
for some bimodule map (called the `generalised braiding') $\sigma:\Omega^1\tens_A\Omega^1\to \Omega^1\tens_A\Omega^1$. This map is not additional data, being determined if it exists by the 2nd Leibniz rule. Thus,  among all connections in the usual sense for a left module, bimodule connections are a nice subset where the right module structure is also respected. The notion goes back to \cite{DV2,Mou}. A connection is {\em metric compatible} if 
\begin{equation}\label{metcomp} \nabla g:=(\nabla\tens\id+ (\sigma\tens \id)(\id\tens\nabla))g=0\end{equation}
using the fact that bimodules equipped with a bimodule connection form a monoidal category. Thus,  $\Omega^1\tens_A\Omega^1$ inherits a tensor product connection given as shown applied to $g$. Finally, a  connection has  torsion and curvature
\begin{equation}\label{torcurv} T_\nabla=\wedge\nabla-\extd:\Omega^1\to \Omega^2,\quad R_\nabla=(\extd\tens\id-\id\wedge\nabla)\nabla:\Omega^1\to \Omega^2\tens_A\Omega^1\end{equation}
and is called {\em quantum Levi-Civita} (QLC) if it is metric compatible and $T_\nabla=0$. There is also a weaker notion
than metric compatibility, namely vanishing of\cite{Ma:rie}
\begin{equation}\label{cotor} coT_\nabla=(\extd\tens\id-\id\wedge\nabla)g\in\Omega^2\tens_A\Omega^1\end{equation}
and a connection which is torsion and cotorsion free is called {\em weak quantum Levi-Civita} (WQLC). One can show that QLC implies WQLC but the latter concept is more general and does not even need the connection to be a bimodule one.

To complete our lightning review of quantum Riemannian geometry, a working definition of the Ricci tensor can be defined relative to a splitting $i:\Omega^2\to \Omega^1\tens_A\Omega^1$ of the wedge product, i.e. a bimodule map such that $\wedge\circ i=\id$. We use this to lift $R_\nabla$ to a map $\Omega^1\to \Omega^1\tens_A\Omega^1\tens_A\Omega^1$ and then take a trace to define\cite{BegMa1} 
\begin{equation}\label{ricci} {\rm Ricci}=( (\ ,\ )\tens\id)(\id\tens (i\tens\id)R_\nabla)g\in \Omega^1\tens_A\Omega^1\end{equation}
and Ricci scalar $S=(\ ,\ ){\rm Ricci}$. These are both $-{1\over 2}$ of their classical values defined by contacting indices of the Riemann tensor,  but our conventions are more natural and work over any field including characteristic 2. It is not clear how to obtain a conserved Einstein tensor in this context, which is why these are working definitions in the absence of a deeper understanding of these matters. 

Also in the general theory, a calculus is {\em inner} if there is an element $\theta\in \Omega^1$ such that $\extd a=[\theta,a]$. The exterior algebra is inner if this also holds in higher degree with the graded-commutator $\extd\omega=[\theta,\omega\}$. If the exterior algebra is inner then all bimodule connections on $\Omega^1$ take the form\cite{Ma:gra} 
\begin{equation}\label{nablatheta} \nabla\omega=\theta\tens\omega-\sigma(\omega\tens\theta)+\alpha(\omega)\end{equation}
for any bimodule maps $\sigma$ mapping as above and $\alpha:\Omega^1\to \Omega^1\tens_A\Omega^1$. This is torsion free iff
\begin{equation}\label{torfreesig} \wedge(\id+\sigma)=0,\quad \wedge\alpha=0\end{equation}
and metric compatible iff
\begin{equation}\label{metcompsig}\theta\tens g-\sigma_{12}\sigma_{23}(g\tens\theta)+(\alpha\tens\id)g+\sigma_{12}(\id\tens\alpha)g=0.\end{equation}

\subsection{Discrete quantum geometry}\label{secX}  There is no problem to specialise the above to the discrete case where $A=k(X)$. We recall this briefly, again from \cite{Ma:gra,BegMa}. Here the different $\Omega^1$ correspond to  different directed graphs with $X$ as vertex set, i.e., to the specification of a set of arrows ${\rm Arr}\subseteq X\times X\setminus {\rm diag}(X)$, where ${\rm diag}(X)$ is the diagonal embedding. Its elements $x\to y$ span $\Omega^1$ as a vector space and the bimodule structure and exterior derivative are
\[ f.x\to y=f(x)x\to y,\quad x\to y.f= f(y)x\to y,\quad \extd f=\sum_{x\to y\in {\rm Arr}}(f(y)-f(x))x\to y\]
for all $f\in A$. If Arr is finite then 
\[ \theta=\sum_{x\to y\in{\rm Arr}}x\to y\]
makes the calculus inner. For $\Omega^1_{uni}$, the graph is the complete one where the ${\rm Arr}=X\times X\setminus{\rm diag}(X)$. Given $\Omega^1$, it follows that $n$-fold tensor products $\Omega^1\tens_A\cdots\tens_A\Omega^1$ are spanned by $n$-step paths,
\[ {\rm Arr}^{(n)}=\{x_0\to x_1\to \cdots\to x_n\ |\ x_i\in X\}.\]
Here ${\rm Arr}^{(1)}={\rm Arr}$ and formally, ${\rm Arr}^{(0)}=X$. The multistep arrow sets partition into subsets  ${}_x{\rm Arr}^{(n)}_y$ where endpoints $x_0=x$ and $x_n=y$ are fixed. It can be shown that a  quantum metric in this context necessarily takes the more specific form \cite{Ma:gra}\cite[Prop 1.28]{BegMa}
\[ g=\sum_{x\to y} g_{x\to y} x\to y\to x\]
for non-zero `arrow weights' $g_{x\to y}$, forced by the requirement of a bimodule inverse, given by $(x\to y, y'\to x'\)=\delta_{x,x'}\delta_{y,y'}\delta_x$.  Note that a quantum metric in this invertible sense requires the directed graph to be bidirected, i.e an undirected graph with arrows understood in both directions. We say that $g$ is {\em edge-symmetric} if $g_{x\to y}=g_{y\to x}$ for all $x\to y$. 

For a given graph, the maximal prolongation exterior algebra is a quotient of the tensor algebra over $A$ on $\Omega^1$  by the sub-bimodule
\[ \CN_{max}=\<\sum_{y}p\to y\to q\ |\ p,q\in X, p\ne q, p\to\kern-8pt/\ q\>\]
but \cite{BegMa} identifies two further canonical quotients of interest, where we quotient by the sub-bimodules
\[ \CN_{max}\subseteq \CN_{med}=\<\sum_{y}p\to y\to q\ |\ p,q\in X, p\ne q\>\subseteq  \CN_{min}=\<\sum_{y}p\to y\to q\ |\ p,q\in X\>.\]
These are generated by the elements shown, i.e. we are imposing relations for all $p,q$ as indicated, with corresponding surjections $\Omega_{max}\to \Omega_{med}\to \Omega_{min}$. If $\Omega^1=\Omega^1_{uni}$ then $\CN_{max}=\{0\}$ as all $p\ne q$ are connected by an arrow, and $\Omega_{max}=\Omega_{uni}$, the universal calculus on $A$. The following is essentially in \cite{BegMa} but now formulated more strongly.  

\begin{lemma}cf \cite{BegMa}. Let {\rm Arr} be finite. Then (i) $\Omega_{med}$ is universal among quotients of $\Omega_{max}$ that remain inner by $\theta$. (ii) $\Omega_{min}$ is universal among quotients of $\Omega_{med}$ that have $\wedge(g)=0$ for the canonical {\em Euclidean metric} $g=\sum_{x\to y}x\to y\to x$.
\end{lemma}
\proof (i) For any unital algebra, if $\Omega^1$ is inner by $\theta$ then 
\[ 0=\extd^2a=\extd(\theta a-a \theta)=(\extd \theta) a- \theta\extd a-(\extd a)\theta-a\extd\theta= [\extd\theta, a]- \theta^2 a+\theta a\theta-\theta a\theta+a\theta^2\]
 requires $\extd \theta-\theta^2$ to be central. If $\Omega^2$ is inner by $\theta$ then $\extd \theta-\theta^2=\theta^2$, so this  needs to be central.  In our graph case, this needs $\theta^2=\sum {x\to y\to z}$ to have nonzero terms in $\Omega^2$ only when $z=x$, which requires (on product by $\delta_p(\ )\delta_q$ for $p\ne q$) precisely the relations of $\Omega_{med}$. (ii) If $g$ is the Euclidean metric then $\wedge(g)=\sum {x\to y\to x}=0$ requires (on multiplication by $\delta_p$ for each $p\in X$) precisely the further relations of $\Omega_{min}$. \endproof

The choice of exterior algebra affects the torsion (and cotorsion) equations. For bimodule connections, the requirement of $\alpha,\sigma$ to be bimodule maps already dictates that
\[ \alpha(x\to y),\ \sigma (x\to w\to y) \in {}_x{\rm Arr}^{(2)}_y\]
for all 1-steps $x\to y$ and all 2-steps $x\to w\to y$ respectively. We then require for zero torsion that 
\[ \alpha(x\to y),\  x\to w\to y+\sigma (x\to w\to y)\in \CN\]
according to the choice of $\CN$ for the degree 2 relations, for all 1-steps and 2-steps. 

Finally, of interest to us is what happens in the special case $k=\F_2$. The first observation is that since $g_{x\to y}\ne 0$ for non-degeneracy of the  quantum metric, we must have $g_{x\to y}=1$ for all $x\to y$, i.e. the above Euclidean metric is the unique quantum metric on a discrete graph calculus. In the inner case, we also have $\extd\theta=2\theta^2=0$ automatically. The bimodule map and zero torsion requirements for $\alpha,\sigma$ now tend to be highly restrictive in the case by case analysis for specific graphs, as we shall in Section~\ref{secex}. It will also be more practical to solve the cotorsion equation to find the WQLCs before imposing full metric compatibility for the moduli of QLCs. 

\section{Boolean algebras and de Morgan duality for differentials}  \label{secboole}  

In principle, most of this section is equivalent to a specialisation of Section~\ref{secX} to the case of $k=\F_2$, with a subset $a\subseteq X$ corresponding to the characteristic function $\chi_a$ which is 1 on the subset $a$ and zero elsewhere. There is, however, a substantial change in language to the level of the power set $P(X)$ of a set $X$, taking  work to disentangle.  It is also natural to use the subset approach when the graph arrow set {\rm Arr} is infinite, which is then slightly more general than the algebraic approach. 

\subsection{Differential Venn diagrams.} \label{secboo}

\begin{figure}
\[  \includegraphics[scale=.8]{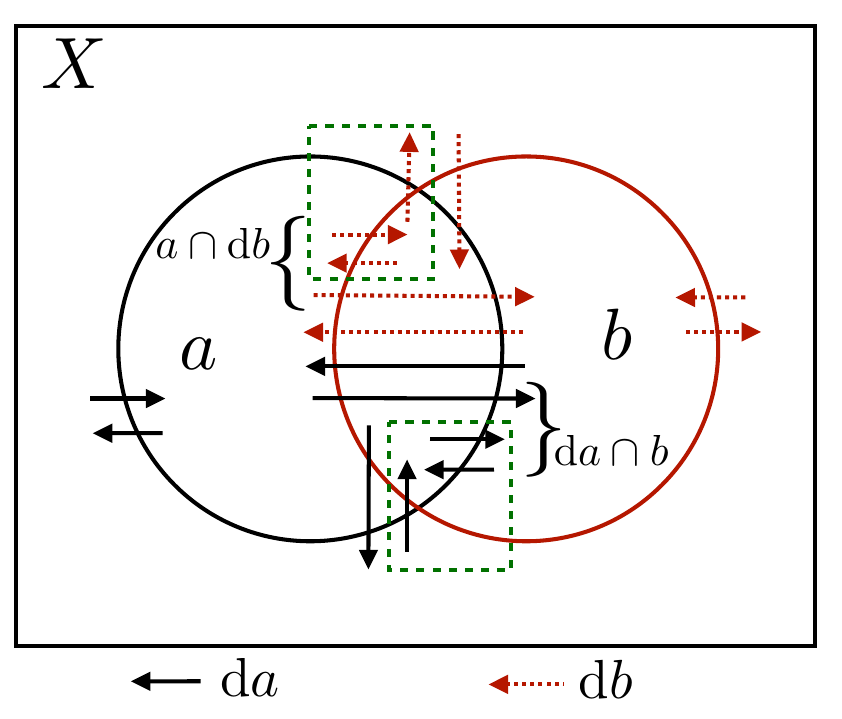}\]
\caption{Venn diagram to check that $\extd(a\cap b)=(\extd a\cap b)\oplus( a\cap \extd b)$. Here $\extd a\cap b$ are arrows in $\extd a$ (solid, black) with tips in $b$ and $a\cap\extd b$ are arrows in $\extd b$ (dotted, red) with tails in $a$. From their union we exclude those in their intersection, leaving the boxed arrows with one end in $a\cap b$ and one end out, i.e. $\extd(a\cap b)$. \label{dab}}
\end{figure}

In the case of functions on a discrete space $X$, the possible $\Omega^1$ are classified by the possible directed graphs with vertex set $X$. So from now on we fix both a set $X$ and a set ${\rm Arr}=\{x\to y\}$ of arrows between some distinct elements of $X$. Here $\Omega^1$ has basis labelled by the arrows and over $\F_2$ each basis element appears or doesn't appear in an element $\omega\in \Omega^1$, so we can identify $\Omega^1(P(X))=P({\rm Arr})$ as the set of subsets of the arrow set of the graph with its $\oplus$ addition law as a Boolean algebra in its own right. Translating the usual finite-difference formulae recalled in Section~\ref{secX} back to $P(X)$, it is easy to see that we  find the following noncommutative bimodule and differential structure:
\begin{gather}\label{setcalc} a\cap\omega:=\{{\rm arrows\ in\ }\omega\ {\rm with\ tail\ in\ }a\},\quad \omega\cap a:=\{{\rm arrows\ in\ }\omega\ {\rm with\ tip\ in\ }a\}\nonumber\\
 \extd a=\{{\rm arrows\ with\ one\ end\ in\ }a\ {\rm and\ other\ end\ in\ }\bar a\},\end{gather}
where $a\subseteq X$ and $\bar a$ is its complement. We extend the usual meaning of $\cap$ to apply between subsets of $X$ and subsets of ${\rm Arr}$ as indicated (but note that this is not commutative) and we use these extensions for the bimodule product, so $a.\omega=a\cap\omega$ and $\omega.a=\omega\cap a$. Thus $\extd a$ is the set of arrows that cross the boundary of $a$ in a Venn diagram.  It's a nice check using Venn diagrams that $\extd$ is indeed a derivation, see Figure~\ref{dab}.  This property in terms of $\cap,\cup$ on $P({\rm Arr})$ is
\[ \extd (a\cap b)=(\extd a\cap b)\oplus (a\cap \extd b)= ((\extd a\cap b)\cup (a\cap \extd b))\cap\overline{(\extd a\cap b)\cap(a\cap \extd b)}\]
where $(\extd a)\cap b)\cap(a\cap \extd b)$ means arrows arrows that cross both $a$ and $b$ boundaries  and have tip in $b$ and tail in $a$, i.e.  the two shown that connect $a\cap\bar b\to \bar a\cap b$ in the figure; we exclude these. 

 We also define $\theta:={\rm Arr}$ as the identity element of $P({\rm Arr})$ and then each subset $a$ partitions the set of all arrows as
 \[ \theta={\rm Arr}=(a\cap\theta \cap a) \oplus \extd a \oplus( \bar a\cap\theta \cap \bar a))\]
into subsets of arrows that, respectively,  lie entirely within $a$ (i.e. the restricted graph on $a\cap X$),  or cross the boundary, or arrows that lie entirely outside $a$ (i.e. the restricted graph on $\bar a\cap X$).  Moreover,
\[ \extd a= (\theta\cap a)\oplus (a\cap \theta)= \theta.a + a.\theta\]
in a more algebraic language for the bimodule products and addition, i.e., the calculus is inner via $\theta$.

We will also need a choice of $\Omega$, for which we take the `maximal prolongation' of $\Omega^1$ (basically, products of 1-forms modulo some minimal set of relations) or further quotients of our choice.  Firstly
\[ \Omega^1\tens_{P(X)}\Omega^1=P({\rm Arr}^{(2)}),\quad \omega\tens_{P(X)} \eta=\{{\rm 2-steps\ starting\ in\ }\omega\ {\rm and\ ending\ in\ }\eta\}\]
where ${\rm Arr}^{(2)}=\{x\to y\to z\}$ denotes the set of 2-step arrows in $X$, i.e. $\tens_{P(X)}$ is the {\em concatenation} of compatible arrows. This is a $P(X)$-bimodule with $a\cap$ and $\cap a$ defined as in (\ref{setcalc}) with `tail' and `tip' now referring to the initial tail or the final tip. Recalling that  ${}_p{\rm Arr}^{(2)}_q$ denotes the set of 2-step arrows between fixed $p,q$, consider the collection of subsets
\begin{equation}\label{N} \CN_{max}=\{ {}_p{\rm Arr}^{(2)}_q| p\to\kern-9pt/\  q\} \subseteq \CN_{med}=\{{}_p{\rm Arr}^{(2)}_q|p\neq q\}\subseteq \CN_{min}=\{{}_p{\rm Arr}^{(2)}_q\},\end{equation}
where the first collection runs over $p\ne q$ for which there is no arrow $p\to q$. We then define $\Omega^2_{max},\Omega^2_{med},\Omega^2_{min}$ by a quotient of $P({\rm Arr}^{(2)})$ by an equivalence relation defined by $\omega\sim\eta$ if $\omega\oplus\eta$ is the union of members of the relevant collection $\CN$ of subsets of ${\rm Arr}^{(2)}$.  One can extend this to all forms but we will need only $\Omega^2$. The max one would be the maximal prolongation in the algebraic setting and the other two are successive quotients, but we take the view that they are defined directly as specified. The latter two are inner with the same $\theta$ as above. Once we have specified the 2-forms, we set
\begin{equation}\label{dom} \extd \omega=\{{\rm 2-steps\ where\ one\ step\ is\ in\ }\omega{\rm \ and\ the\ other\ step\ is\ not}\}\end{equation}
but with the output viewed up to the chosen equivalence. One can check for example that 
\[ \extd\extd a\sim \emptyset\]
for all $a\subseteq X$. Here the left hand side consists of all 2-steps where one step crosses the boundary of $a$ and the other does not cross the boundary of $a$. If we fix $p\in a$ and $q\in\bar a$, for example and if there is such a 2-step $p\to x\to q$ then all $x$ meet the criterion so all of ${}_p{\rm Arr}^{(2)}_q$ is included. Similarly for $p\in \bar a$ and $q\in a$. 

We have outlined our constructions of our three exterior algebras directly in the set-theoretic setting. If $X$ and hence ${\rm Arr}$ for a connected calculus are infinite then we can impose a further surjectivity condition to fit with the usual algebraic setting, to the effect that every subset of {\rm Arr} is the $\oplus$ of a finite number of 1-forms of the form $a\extd b$ for $a,b\subseteq X$. It is also possible, and could be more natural here, to proceed without such a surjectivity condition\cite{MaTao2}. In practice, this will not be an issue as our sets will be finite. . 

\subsection{De Morgan duality for differential forms} \label{secDM}

The classical de Morgan's theorem is that for any equality in a Boolean algebra we can swap 
\[ a \leftrightarrow \bar a,\quad \cap \leftrightarrow \cup,\quad 1=X \leftrightarrow  0=\emptyset\]
and still have a valid equality. Moreover, $\overline{a\cap b}=\overline a\cup \overline b$ holds for all $a,b$.  In this section, we want to see how this duality extends to differential forms. 

The first thing to note is that complementation does not respect addition by $\oplus$ on $P(X)$ so it cannot be expressed as any kind of operator on this as a vector space over $\F_2$. Rather, we define $\bar P(X)$ as again the power set of subsets of $X$ but now with product given by $\cup$ and addition given by 
the de Morgan dual-exclusive-OR  (built using $\cap,\cup$ swapped), namely what we call  {\em inclusive AND},
\[ a\bar\cdot b:= a\cup b,\quad a\bar\oplus b:=(a\cap b)\cup\overline{a\cup b}=(a\cap b)\cup (\bar a\cap \bar b)=\overline{(a\cup b)\cap\overline{a\cap b}}=\overline{a\oplus b}.\]
One can check that this again makes the power set of $X$ into an algebra over $\F_2$ (as it must by de Morgan's theorem) and that we now have an isomorphism of algebras
\[ \bar{\ }:P(X)\to \bar P(X).\]
Here $\bar a\oplus\bar b= a\oplus b$ so that $\overline{\bar a\oplus\bar b}=\overline{a\oplus b}=a\bar\oplus b$ as required for linearity over $\F_2$. 

\begin{figure}
\[ \includegraphics[scale=.8]{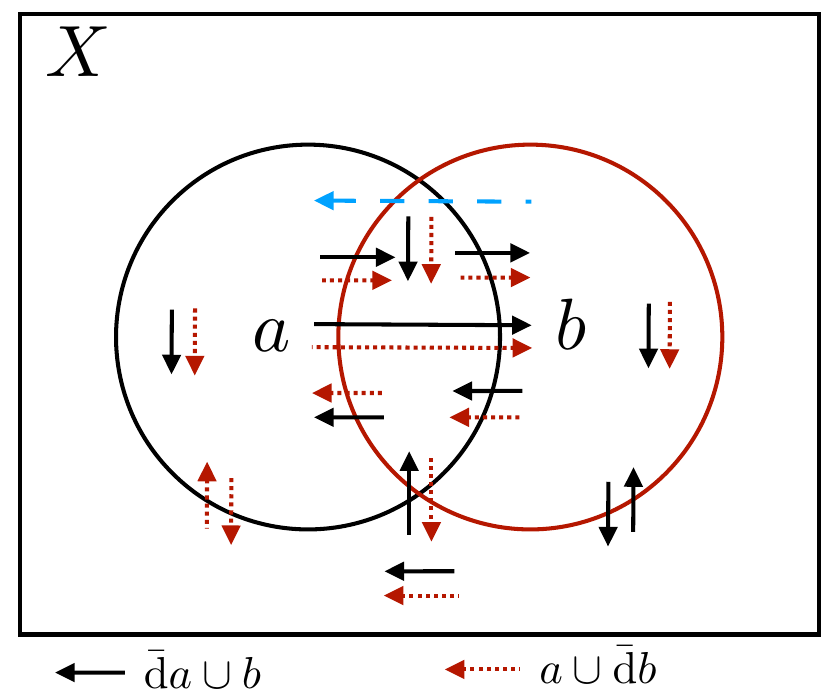} \]
\caption{Venn diagram to check that $\bar\extd(a\cup b)=(\bar\extd a\cup b)\bar\oplus ( a\cup \bar\extd b)$. Here  $\bar\extd a\cup b$ (solid, black) are arrows wholly in or out of $a$ or have tips in $b$, while $a\cup \bar\extd b$ (dotted, red) are arrows wholly in or out of $b$ or have tails in $a$. The parallel arrows are in both  subsets and we see that they are almost all the arrows wholly in or out of $a\cup b$, i.e. almost $\bar\extd(a\cup b)$. We need to add the missing type of arrow shown dashed, blue, which is not in the union of all the other arrows shown.  \label{bardab}}
\end{figure}

Next, define the 1-forms  $\bar\Omega^1:=\Omega^1(\bar P(X)):=\bar P({\rm Arr})$ meaning its addition law is by $\bar\oplus$ of subsets of arrows, with bimodule structure and exterior derivative 
\begin{gather}\label{dualcalc}  a\cup \omega=\{{\rm arrows\ in\ }\omega\ {\rm or\ with\ tail\ in\ }a\},\quad \omega\cup a=\{{\rm arrows\ in\ }\omega\ {\rm or\ with\ tip\ in\ }a\}\nonumber\\
\bar \extd a= \{{\rm arrows\ wholly\ in\  }a\ {\rm or\ wholly\ in\ }\bar a\}=\overline{\extd a}\end{gather}
where complementation of a subset of arrows is in {\rm Arr}. We extended $\cup$ to apply between subsets of $X$ and subsets of ${\rm Arr}$ as stated and a little thought shows that
\begin{equation}\label{extdemorg} \overline{a\cup\omega}=\bar a\cap\bar\omega,\quad \overline{\omega\cup a}=\bar\omega\cap \bar a.\end{equation}
We use this extended $\cup$ for the bimodule structure of $\bar P(X)$, so $a\bar\cdot \omega=a\cup\omega$ and $\omega\bar\cdot a=\omega\cup a$. Figure~\ref{bardab} checks that this indeed obeys the derivation rule for a differential calculus on $\bar P(X)$. However, this must be the case by de Morgan's theorem in view of the symmetry between $\cup$ and $\cap$. In terms of $\cup,\cap$ this is
\[ \bar\extd(a\cup b)=(\bar\extd a\cup b)\bar\oplus (a\cup \bar\extd b)=  ((\bar\extd a\cup b)\cap (a\cup \bar\extd b))\cup\overline{(\bar\extd a\cup b)\cup (a\cup\bar\extd b)}.\]
Here $(\bar\extd a\cup b)\cup (a\cup\bar\extd b)$ means arrows wholly on our out of $a$ or with tip in $b$ or wholly in or out of $b$ or with tail in $a$. 

One can also check that this calculus is inner with $\bar\theta=\emptyset$ of arrows. Thus
\begin{align*}(\emptyset\cup a)\bar\oplus(a\cup\emptyset)&=\{{\rm arrows\ with\ tip\ in\ }a\}\bar\oplus\{{\rm arrows\ with\ tail\ in\ }a\}\\
&=\{{\rm arrows\ wholly\ in\ }a{\rm\ or\ }\bar a\}=\bar\extd a\end{align*}
using the above definition of $\bar\oplus$. 

\begin{proposition}\label{barOmega1} Let $X$ be a graph. The algebra isomorphism $\bar{\ }:P(X)\to \bar P(X)$ with their respective differential structures is a diffeomorphism.
\end{proposition}
\proof The key observation is that $\extd,\bar\extd$ as defined are symmetric between $a,\bar a$, so in particular we have  $\bar\extd \bar a=\bar\extd a=\overline{\extd a}$ so that complementation of arrows forms a commutative diagram 
\[ \begin{array}{ccc}P({\rm Arr})  &{\buildrel \bar{\ }\over \longrightarrow} & \bar P({\rm Arr})   \\ \extd\uparrow & & \uparrow \bar\extd\\ P(X) & {\buildrel \bar{\ }\over \longrightarrow}& \bar P(X).\end{array}\]
The top map is a bimodule map in the sense $\overline{a.\omega}=\bar a\, \bar\cdot\, \bar\omega$ and similarly on the other side, by the observation (\ref{extdemorg}) already given. 
 \endproof

Next, whereas $\omega\tens_{P(X)}\eta$ is the set of possible concatenations or a kind of intersection of a tip in $\omega$ and a tail in $\eta$, we define the dual {\em coconcatenation} of subsets of arrows
\[ \omega\tens_{\bar P(X)}\eta=\{{\rm 2-steps\ starting\ in\ }\omega\ {\rm or\ ending\ in\ }\eta\}=\overline{\bar\omega\tens_{P(X)}\bar\eta}\]
and one can check that 
\[ (\omega\cup a)\tens_{\bar P(X)}\eta=\omega\tens_{\bar P(X)}(a\cup\eta)\]
as both sides are arrows that start in $\omega$ or end in $\eta$ or have middle vertex in $a$. For the addition law we use $\bar\oplus$ and we define a bimodule structure on 2-steps by extending $a\cup$ and $\cup a$ in (\ref{dualcalc}) to 2-steps with `tail' and `tip' referring to the initial tail or the final tip. We still have (\ref{extdemorg}) with this extension. In this way we  identify 
\[ \bar \Omega^1\tens_{\bar P(X)}\bar \Omega^1= \bar P({\rm Arr}^{(2)})\]
 and by construction the complementation map
\[  \bar{\ }: P({\rm Arr}^{(2)})=\Omega^1\tens_{P(X)}\Omega^1\to \bar P({\rm Arr}^{(2)})=\bar\Omega^1\tens_{\bar P(X)}\bar\Omega^1\]
intertwines the bimodule structures in same way as for $\bar{\ }:\Omega^1\to \bar\Omega^1$ in Proposition~\ref{barOmega1}, namely  $\overline{a\cap\omega\tens_{P(X)}\eta}=\bar a\cup\overline{\omega\tens_{P(X)}\eta}$ and similarly on the other side.  

\begin{lemma} $\bar{\ }$ descends to the relevant max, med, min prolongations in a way that commutes with $\extd$, $\bar\extd$ on degree 1. 

\end{lemma} 
\proof We let $\CN$ be one of the collections (\ref{N}). Its elements are the $\oplus$ of any subset $Y$ of the allowed $(p,q)$ in the relevant collection and such an element maps to  $\overline{\oplus_{(p,q)\in Y}{}_p{\rm Arr}^{\rm (2)}_q}=\bar\oplus_{(p,q)\in Y}{}_p{\rm Arr}^{\rm (2)}_q$ as the corresponding element of $\bar\CN$. The latter is defined by the same collections as $\CN$ but with elements constructed in $\bar P({\rm Arr}^{(2)})$ using $\bar\oplus$. Then by construction $\bar{\ }$ descends to $\Omega^2$ where we quotient by $\CN$ on mapping to  $\bar\Omega^2$ where we quotient by $\bar\CN$.  That the differentiability diagram for $\extd,\bar\extd$ commutes follows similarly to the proof for $P(X)$ given the form of $\extd$ in (\ref{dom}) and the dual version
\begin{equation}\label{bardom} \bar\extd \omega=\{{\rm 2-steps\ wholly\ in\ }\omega{\rm \ or\ wholly\ out}\}.\end{equation}\endproof

\subsection{Elements of quantum Riemannian geometry on $P(X)$} 

We continue in the case of $X$ a graph with $A=P(X)$ and $\Omega^1=P({\rm Arr})$. Now we suppose the graph is bidirected (so for every arrow $x\to y$ there is an arrow $y\to x$). Then the unique quantum metric is
\begin{equation}\label{gset} g=\{x\to y\to x\ |\ x\to y\in {\rm Arr}\}=\cup_p\, {}_p{\rm Arr}^{(2)}_p\in \Omega^1\tens_{P(X)}\Omega^1\end{equation}
 i.e. all 2-steps that go to another point and come back. 
 
A connection has to map subsets of ${\rm Arr}$ to subsets of ${\rm Arr}^{(2)}$ subject to certain properties and with constructions transferred from the algebraic theory. In particular, at least if $X$ is finite,  bimodule connections will be determined by bimodule maps  $\alpha:P({\rm Arr})\to P({\rm Arr}^{(2)})$ and $\sigma:P({\rm Arr}^{(2)})\to P({\rm Arr}^{(2)})$ as (\ref{nablatheta}), with $\tens_{P(X)}$ now understood as in Section~\ref{secboo}. To be torsion free now appears as $\alpha,\sigma+\id$ having image in the chosen $\CN$. Thus, for $\CN_{min}$, each element in the image should be the union of subsets of the form ${}_p{\rm Arr}^{(2)}_q$. Likewise, we can use (\ref{metcompsig}) for metric compatibility or proceed directly from (\ref{metcomp}), adapted to the subset form.

\section{Small connected Boolean Riemannian geometries}\label{secex}

For $X$ a singleton, the only graph is no edges and $\Omega^1=0$, so there is no quantum metric and no Riemannian geometry. In  Sections~\ref{sec2pt}-\ref{secpoly}, we fully solve the $|X|=2,3$ cases, the square case for $|X|=4$ and $n$-gon case for $n>5$. We number the vertices $0,1,2,\cdots$ and adopt shorthand notation $01=0\to 1$ and $010=0\to 1\to 0$, etc. for the arrows and multistep arrows. We also use these as representatives in quotient spaces or we may write $[010]$ etc., if we explicitly want to indicate an  equivalence class. We proceed in the $\F_2$-algebra setting but exhibit the conversion to Boolean subset form, collected for the polygon case in  Section~\ref{secsub}.  Section~\ref{secdem} then exhibits an example of de Morgan duality at the geometric level. 

\subsection{All quantum geometries on 2 points}\label{sec2pt} For $X=\{0,1\}$ the only connected graph is $0-1$ up to relabelling and hence $\Omega^1=\<01,10\>=\Omega^1_{uni}$ since the latter is 1-dimensional over the algebra.  We have $\CN_{max}=\CN_{med}=\{0\}$ so that $\Omega^2_{max}=\Omega^2_{med}=\Omega^2_{uni}=\Omega^1\tens_A\Omega^1=\<010,101\>$. 

The unique metric over $\F_2$ is therefore
\[ g=010+101=\theta\tens\theta;\quad \theta=01+10\]
and is not quantum symmetric for $\Omega_{uni}$.  We also have $\Omega^2_{min}=\{0\}$, for which $\theta^2=0$ and $\wedge(g)=0$, but we do not assume that we work with this. 

\begin{lemma} There is a unique bimodule connection on $0-1$ with invertible $\sigma$, 
\[ \nabla(01)=\nabla(10)=g,\quad \sigma=\id.\]
This is also the unique metric compatible bimodule connection and has $T_\nabla=0$ (so a QLC) and $R_\nabla=0$.
\end{lemma}
\proof The calculus is inner with $\theta=01+10$, so connections are given by $\alpha,\sigma$. We must have $\alpha(01)=\alpha(10)=0$ for $\alpha$ to be a bimodule map. Similarly, for a bimodule map $\sigma$, we must have $\sigma(010)=a 010$, $\sigma(101)=b101$ for coefficients in the field. Then $\nabla(01)=101+a010$ and $\nabla(10)=010+b101$, so that
\begin{align*}\nabla g&=\nabla(01)10+\nabla(10)01+(\sigma\tens\id)(01\nabla(10)+10\nabla(01))\\
&=1010+0101+(\sigma\tens\id)(b0101+a1010)=1010+0101+ab0101+ab1010\end{align*}
which vanishes if and only if $a=1,b=1$, i.e. $\sigma=\id$. This is also the only invertible $\sigma$. One can then check the torsion and curvature, e.g.,
\[ T_\nabla(01)=(\wedge\nabla+\extd)(01)=\wedge(g)+(01+10)01+01(01+10)=0\]
\[ R_\nabla(01)=(\extd\tens\id-\id\wedge\nabla)\nabla(01)=coT_\nabla(01)=0\]
since QLC implies WQLC (one can also directly compute that this vanishes).  \endproof

In terms of subsets of the relevant arrow spaces, we have  $\Omega^1=P({\rm Arr})=P(\{01,10\})$ and in here $\theta=\{01,10\}={\rm Arr}$. We also have $\CN_{max}=\CN_{med}=\emptyset$ and  $\Omega^2_{uni}=P({\rm Arr}^{(2)})=\Omega^1\tens_A\Omega^1=P(\{010,101\})$, and in here we have
\[ g=\{010,101\}\subset {\rm Arr^{(2)}},\quad \nabla {\rm Arr}=\emptyset,\quad \nabla\{01\}=\nabla\{ 10\}=g.\]
as well as $\sigma=\id$ on ${\rm Arr}^{(2)}$. The content is the same as in terms of linear algebra over $\F_2$, but now in an unfamiliar subset point of view.

\subsection{\bf All quantum geometries on 3 points.}\label{sec3pt} 

For $X=\{0,1,2\}$, there are two connected graphs up to relabelling, the line segment $0-1-2$ and the triangle. 

\subsubsection{Line segment graph $0-1-2$} Here, $\Omega^1=\<01,10,12,21\>$ while 
\[ \CN_{max}=\CN_{med}=\<012,210\>,\quad \CN_{min}=\<012,210,010,101+121,212\>\]
 so that $\Omega^1$ and $\Omega^2_{max}=\Omega^2_{med}$ are 4-dimensional and $\Omega^2_{min}$ is $1$-dimensional as vector spaces. This means that these are all strictly projective and not free modules. We have that $\Omega^1\tens_A\Omega^1=\<010,012,101,121,210,212\>$ is 6-dimensional.  The unique quantum metric is 
\[ g=010+101+121+212\]
and we have $\theta=01+10+12+21$ with $\theta\tens\theta=g+012+210$. This means that $\wedge(g)=\theta^2$ in all cases and vanishes for $\Omega_{min}$. 

\begin{proposition}\label{proplinenabla} There is no metric compatible bimodule connection $\nabla$ on $0-1-2$. There is unique bimodule WQLC for $\Omega_{max}$ and four for $\Omega_{min}$, namely
\begin{gather*} \nabla(01)=101+010,\quad \nabla(10)=010+210+\alpha101+(1+\alpha)121,\\
 \nabla(21)=121+212,\quad \nabla(12)=212+012+(1+\delta)101+\delta121,\\\sigma(010)=010,\quad\sigma(012)=\sigma(210)=0,\quad\sigma(212)=212,\\
\sigma(101)=\alpha 101+(1+\alpha)121,\quad \sigma(121)=(1+\delta)101+\delta 121\end{gather*}
for $\alpha,\delta\in\F_2$ (here $\alpha=\delta=1$ is also the $\Omega_{max}$ solution).    
\end{proposition}
\proof We only have $\alpha(01)=\alpha(10)=\alpha(12)=\alpha(21)=0$ for a bimodule map, while the form of a bimodule map $\sigma$ is necessarily 
\begin{gather*} \sigma(010)=a 010,\quad \sigma(012)=b 012,\quad \sigma(210)=c210,\quad \sigma(212)=d212,\\ 
\sigma(101)=\alpha 101+\beta 121,\quad \sigma(121)=\gamma 101 +\delta 121\end{gather*}
for some coefficients in field. Then
\begin{gather*} \nabla(01)=101+a010+b012,\quad \nabla(10)=010+210+\alpha101+\beta121,\\ 
 \nabla(21)=121+c210+d212,\quad \nabla(12)=012+212+\gamma101+\delta121.\end{gather*}
We then compute
\begin{align*}\nabla g&=\nabla(010+101+121+212)\\
&=\nabla(01)10+\nabla(10)01+\nabla(12)21+\nabla(21)12+(\sigma\tens\id)(01\nabla(10)+10\nabla(01)+21\nabla(12))\\
&=1010+2101+0101+0121+2121+1212\\
&\quad+(\sigma\tens\id)(\alpha0101+\beta0121+a1010+b1012+d1212+c1210+\gamma2101+\delta2121)\\
&=(1+\alpha a+\gamma c)1010+(1+\alpha a)0101+(1+\beta b)0121+(1+\beta b+\delta d)1212\\
&\quad+(1+\delta d)2121+(1+\gamma c)2101+(\gamma d+\alpha b)1012+(\beta a+\delta c)1210\end{align*}
as an element of $\Omega^1\tens_A\Omega^1\tens_A\Omega^1$ (we suppressed $\tens$ in displaying the elements). We used $\nabla$ and then $\sigma$. Vanishing of all 6 terms needs on the one hand $\gamma c=0$ and on the other hand $\gamma c=1$, so is not possible. So there are no metric compatible bimodule connections.

 If we project to $\Omega^2_{max}\tens_A\Omega^1$, however, we kill 012 and 210 in the first factor, so we just need
\[ \alpha a=\delta d=1,\quad \gamma c=\beta b=0,\quad \gamma d=\alpha b,\quad \beta a=\delta c.\]
This needs $\alpha=\delta=a=d=1$ and either (a) $\beta=\gamma=b=c=0$ (b) $c=\beta=1$, $b=\gamma=0$ or (c) $b=\gamma=1$, $b=c=0$. None of these have $\sigma$ invertible. If we project further to $\Omega_{min}$ then we just need
\[ (\alpha +\beta)a +(\gamma+\delta)c=(\alpha+\beta)b+(\gamma+\delta d)=1\]
with many more solutions. As an aside, if we restrict to invertible $\sigma$ then $a=b=c=d=1$ and we just have the two equations $\alpha+\beta+\gamma+\delta=\alpha\delta+\beta\gamma=1$ (the latter for invertibility), which has four solutions. So there are four cotorsion free connections with $\sigma$ invertible and many more otherwise, for $\Omega_{min}$.

For the torsion, we compute $\extd =[\theta,\ ]$ in $\Omega_{max}$ as 
\[ \extd(01)=\extd(10)=010+101,\quad \extd(12)=\extd(21)=121+212\]
and comparing with $\wedge\nabla$, we see that $T_\nabla=0$ if and only if $a=d=\alpha=\delta=1$ and $\beta=\delta=0$ which intersects with solution (a) in the preceding paragraph for the cotorsion to vanish. This is the unique bimodule WQLC here. For $\Omega_{min}$ we need only $a=d=\alpha+\beta=\gamma+\delta=1$. The cotorsion vanishing then needs $b=c=0$, giving the four bimodule WQLCs in this case. \endproof

In fact, all four connections are curved and we will compute their Ricci tensors. This requires the metric inner product  \[ (01,10)=\delta_0,\quad (10,01)=(12,21)=\delta_1,\quad (21,12)=\delta_2\]
in terms of Kronecker $\delta$-functions in $A$, and zero for the remaining combinations from the basis. Next, a bimodule map $i:\Omega^2_{max}\to \Omega^1\tens_A\Omega^1$ such that $\wedge\circ i=\id$, must have the form
\[ i([010])=010,\quad i([212])=212,\quad i([101])=101,\quad i([121])=121\]
understood as tensor products on the right. This lead to a uniquely-defined Ricci tensor. The same requirements for $i:\Omega^2_{min}\to \Omega^1\tens_A \Omega^1$ has two solutions,
\[ i_1([101])=101,\quad i_2([101])=121,\]
leading to two  Ricci tensors for each connection with this exterior algebra, according to the choice of lift. 

\begin{proposition} On $0-1-2$, the unique WQLC for $\Omega_{max}$  in Proposition~\ref{proplinenabla} has curvature
\[R_\nabla(01)=R_\nabla(21)=0,\quad R_\nabla(10)=1210,\quad  R_\nabla(12)=1012\]
as elements of $\Omega^2_{max}\tens_A\Omega^1$, and is Ricci flat. The four WQLC's for $\Omega_{min}$ have curvature
\[ R_\nabla(01)=R_\nabla(21)=0,\quad R_\nabla(10)=\alpha1010+(1+\alpha)1012,\quad R_\nabla(12)=(1+\delta)1010+\delta1012.\]
as elements of  $\Omega^2_{min}\tens_A\Omega^1$ with two Ricci tensors in $\Omega^1\tens_A\Omega^1$ and corresponding scalars in $A$,
\[ {\rm Ricci}_1=\alpha010+(1+\alpha)012,\quad S_1=\alpha,\quad {\rm Ricci}_2=\delta212+(1+\delta)210,\quad S_2=\delta\]
according to the lifts $i_1,i_2$ respectively. 
\end{proposition}
\proof  The curvature is 
\begin{align*} R_\nabla(01)&=(\extd\tens\id+\id\wedge\nabla)(101+010)=\extd(01)01+\extd(01)10+10\nabla(01)+01\nabla(10)\\
R_\nabla(10)&=(\extd\tens\id+\id\wedge\nabla)(010+210+\alpha101+(1+\alpha)121)\\
&=\extd(01)10+\extd(21)10+\alpha\extd(10)01+(1+\alpha)\extd(12)21\\
&\quad +01\nabla(10)+21\nabla(10)+\alpha10\nabla(01)+(1+\alpha)12\nabla(21)\\
\end{align*}
and similarly, by $0\leftrightarrow2$ symmetry for the other half. This gives
\begin{gather*}R_\nabla(01)=(1+\alpha)0101,\quad R_\nabla(10)=1210+(1+\alpha)(1010+1212),\\
R_\nabla(21)=(1+\delta)2121,\quad R_\nabla(12)=1012+(1+\delta)(1010+1212)\end{gather*}
of which $\alpha=\delta=1$ is the stated result for the WQLC for this calculus. We also quotient these expressions to obtain the  curvatures for $\Omega_{min}$ as stated.  

For Ricci, we use $i$ to lift the first tensor factor in the output of $R_\nabla$ then use the metric and inverse metric to make a trace. In the first case, 
\begin{align*}{\rm Ricci}&=((01,\  )\tens\id)(i\tens\id)R_\nabla(10)+ ((21,\  )\tens\id)(i\tens\id)R_\nabla(12)\\
&=((01,\  )\tens\id)1210+ ((21,\  )\tens\id)1012,\end{align*}
where we apply $(i\tens\id)R_\nabla$ to the second tensor factor of each term in $g$ and contract the first factor of its output with the corresponding left factor of that term in $g$. In our notation, $i$ turns $1210$ understood as representing $[121]\tens 10$ to $1210$ understood as $12\tens 21\tens10$, for example. We then evaluate the inverse metric to obtain the result stated.  For  the $\Omega_{min}$ WQLCs, we similarly have
\begin{align*} {\rm Ricci}&=((01,\ )\tens\id)(\alpha i([101])10+(1+\alpha)i([101])12)\\
&\quad + ((21,\ )\tens\id)((1+\delta)i([101])10+\delta i([101])12)\end{align*}
which, for $i=i_1$ and $i=i_2$, gives the two Ricci curvatures stated. \endproof

In physics, a non-flat but Ricci flat metric, as here for $\Omega_{max}$, would be a vacuum solution of Einstein's equation (such as a black hole). The Ricci tensors for $\Omega_{min}$ show the geometric meaning of the parameters $\alpha,\delta$. 
 
Finally, in the subset language of Section~\ref{secboole},  $\Omega^1=P(\{01,10,12,21\})$ and in here the inner generator is $\theta=\{01,10,12,21\}$. We have $\Omega^1\tens_A\Omega^1=P(\{010,012,101,121,210,212\})$ and in here the unique quantum metric and the values of the QLCs and the Ricci tensors for $\Omega_{min}$ are
\begin{gather*} g=\{010,101,121,212\},\quad \nabla\{01\}=\{101,010\},\quad\nabla\{21\}=\{121,212\}\\ 
 \nabla\{10\}=\begin{cases}\{010,210,101\}&\alpha=1\\ \{010,210,121\}&\alpha=0\end{cases},\quad \nabla\{12\}=\begin{cases}\{212,012,121\}&\delta=1\\ \{212,012,101\}&\delta=0\end{cases}\\ 
 {\rm Ricci}_1=\begin{cases}\{010\}&\alpha=1\\ \{012\}&\alpha=0\end{cases},\quad {\rm Ricci}_2=\begin{cases}\{212\}&\delta=1\\ \{210\}&\delta=0.\end{cases}\end{gather*}
Here, the exterior algebras that we considered were defined by 
\[ \CN_{max}=\CN_{med}=P(\{012,210\}),\quad \CN_{min}=P(\{012,210,010,101+121,212\}).\]

\subsubsection{Triangle graph}\label{sectri} For  the triangle graph, $\Omega^1=\<01,10,12,21,20,02\>=\Omega^1_{uni}$ since the latter is a free module and 2-dimensional over the algebra. We have
\[ \CN_{max}=\{0\},\quad \CN_{med}=\<021,012,120,102,210,201\>,\]
\[\CN_{min}=\<021,012,120,102,210,201, 010+020, 101+121, 202+212\>\]
so that $\Omega^2_{max}=\Omega^2_{uni}$ is 4-dimensional over the algebra, $\Omega^2_{med}$ is 2-dimensional over the algebra while $\Omega^2_{min}$ is 1-dimensional over the algebra. Here 
\[ \Omega^1\tens_A\Omega^1=\<020,010,012,021,101,121,120, 102, 212, 202,201,210\>\]
 is 4-dimensional over the algebra. The inner generator and the unique quantum metric are
\[ \theta=01+10+12+21+20+02,\quad g=010+020+101+121+202+212.\]
Here $\wedge(g)=0$ for $\Omega_{min}$ but not for the others. 

In fact, the graph here is a Cayley graph for $\Z_3$ and $\Omega^1$ has a natural left-invariant basis $e^\pm$,  with the element $\theta$ and the unique quantum metric in these terms,
\[ e^+=01+12+20,\quad e^-=10+ 21+ 02,\quad  \theta=e^++e^-,\quad g=e^+\tens e^-+e^-\tens e^+.\]
Then $\Omega_{max}$ is generated over $A$ by $e^\pm$ with no relations between them, $\Omega_{med}$ adds the two relations $(e^\pm)^2=0$ while $\Omega_{min}=\Omega(\Z_3)$, the canonical Cayley graph exterior algebra\cite{BegMa} has one further relation $e^+\wedge e^-+e^-\wedge e^+=0$. In all cases, $\extd e^\pm=0$. Here $\Omega(\Z_3)$  is more like a  parallelisable 2-manifold from the point of view of the DGA,  with a top or `volume' form
\[ {\rm Vol}=e^+\wedge e^-=[010+121+202].\]
We will also need the inverse metric given by $(e^\pm,e^\mp)=1$ as an element of $A$ and $(e^\pm,e^\pm)=0$, and a lift $i:\Omega^2(\Z_3)\to \Omega^1\tens_A\Omega^1$, for which there are several choices but two natural (left invariant) ones,
\begin{equation}\label{ipm} i_+({\rm Vol})= e^+\tens e^-,\quad i_-({\rm Vol})= e^-\tens e^+.\end{equation} 
This gives two natural Ricci tensors ${\rm Ricci}_\pm$ according to the choice of lift. By the classification results in \cite{MaPac} (for the $n=3$ algebra B there)  there are in fact four QLCs for $\Omega(\Z_3)$, which we now write much more simply as follows.

\begin{proposition}\label{proptri} On the triangle, for $\Omega_{min}=\Omega(\Z_3)$, the unique metric $g$ has exactly four QLCs, 
\[ \nabla e^+=\alpha e^-\tens e^-,\quad \nabla e^-=\beta e^+\tens e^+,\quad \alpha,\beta\in\{0,1\)\]
with curvature $R_\nabla e^\pm= \alpha\beta {\rm Vol}\tens e^\pm$ and Ricci tensors and Ricci scalar
\[ {\rm Ricci}_+ = \alpha\beta e^-\tens e^+,\quad {\rm Ricci}_-=\alpha\beta e^+\tens e^-,\quad S_\pm=\alpha\beta.\]
\end{proposition}
\proof We let  $\sigma={\rm flip}$ on tensor products of $e^\pm$, which means that $\nabla e^\pm=\alpha(e^\pm)$ with the values as stated. These are clearly bimodule maps when we use the relations $e^\pm f=(R_\pm f)e^\pm$ with $R_\pm(f)(i)=f(i\pm 1)$ for $i=0,1,2$ mod $3$. The connections are torsion free as $\extd e^\pm=0$ and $(e^\pm)^2=0$. For metric compatibility,
\begin{align*}
\nabla(e^+\tens e^-+e^-\tens e^+)&=\alpha e^-\tens e^-\tens e^-+\beta e^+\tens e^+\tens e^++(\sigma\tens\id)(e^+\tens\beta e^+\tens e^++e^-\tens\alpha e^-\tens e^-)\end{align*}
which vanishes since $\sigma$ acts as the identity.  It then follows from the computer results in \cite{MaPac} that these are all the QLCs, since we found 4. In fact, three of them are the $\F_2$ case of the generic field solutions for the triangle in \cite{BegMa} but the curved one where $\alpha=\beta=1$ is specific to $\F_2$. The computations for the curvature and Ricci tensors are straightforward given the simple form of the metric and connections in the $e^\pm$ basis. \endproof

Note that classically, we would lift the volume form  ${\rm Vol}$ to an antisymmetric cotensor with $i=(i_+-i_2)/2$, but we do not have that luxury over $\F_2$. A new approach which we propose in our case is to instead define  `twice the Ricci tensor' ${}^2{\rm Ricci}$ using $i=i_++ i_-$ without the $1/2$ factor (and without the sign as we work over $\F_2$) and the corresponding `twice Ricci scalar' ${}^2S$. This gives
\[ {}^2{\rm Ricci}={\rm Ricci}_++{\rm Ricci}_-=\alpha\beta g,\quad {}^2S=\alpha\beta(\ ,\ )(g)=0.\]
The Ricci expression says that the Boolean algebra with this calculus and QLC  is `Einstein'. 

 Another question concerns the Einstein tensor and a tentative proposal in \cite{MaPac} for quantum geometry over $\F_2$ is ${\rm Eins}:={\rm Ricci}+S g$, i.e. without the $-{1/2}$ which would normally be needed in the 2nd term. 
It was found in \cite{MaPac} for this model that among all lift maps $i$, there are exactly two lifts for which Eins is conserved in the sense $\nabla\cdot{\rm Eins}:=((\ ,)\tens\id)\nabla {\rm Eins}=0$ and we understand these now as precisely $i_\pm$ in (\ref{ipm}) in our $e^\pm$ description. For these lifts,
\[  {\rm Eins}_\pm=\alpha\beta e^\pm\tens e^\mp=\alpha\beta i_\pm({\rm Vol})\]
and 
\[ \nabla ({\rm Eins}_+)=\alpha\beta\nabla (e^+\tens e^-)=\alpha\beta(\nabla e^+\tens e^-+\sigma_{12}(e^+\tens \nabla e^-))=\alpha\beta(e^-\tens e^-\tens e^-+e^+\tens e^+\tens e^+),\]
which  vanishes when we contract the first two factors with the inverse metric $(\ ,\ )$. Similarly for conservation of ${\rm Eins}_-$. Such conservation is also obviously true for their sum ${}^2{\rm Eins}={}^2{\rm Ricci}=\alpha\beta g$. 
Finally, for completeness, we also look at the more general picture of WQLCs.

\begin{proposition}\label{propwqlctri} $\Omega(\Z_3)$ has a 4-functional parameter space of WQLCs,
\begin{gather*} \nabla e^+=\alpha e^-\tens e^-+\gamma e^+\tens e^++\delta g,\quad \nabla e^-=\beta e^+\tens e^++\delta e^-\tens e^-+ \gamma g,\\
\sigma(e^+\tens e^+)=(1+\gamma)e^+\tens e^+,\quad \sigma(e^+\tens e^-)=(1+\delta)e^-\tens e^++\delta e^+\tens e^-,\\
\sigma(e^-\tens e^-)=(1+\delta)e^-\tens e^-,\quad \sigma(e^-\tens e^+)=(1+\gamma)e^+\tens e^-+\gamma e^-\tens e^+\end{gather*}
for $\alpha,\beta,\gamma,\delta\in A$, with  curvature
\begin{gather*} R_\nabla e^+= {\rm Vol}\tens((\alpha\beta+\gamma\delta+\del_-\gamma)e^++\del_+\alpha e^-+\extd \delta),\\ R_\nabla e^-= {\rm Vol}\tens((\alpha\beta+\gamma\delta+\del_+\delta) e^- +\del_-\beta  e^++\extd \gamma).\end{gather*}
The case $\gamma=\delta=0$ and $\alpha,\beta$ constant recovers the QLC's in Proposition~\ref{proptri}.
\end{proposition}
\proof  The possible bimodule maps $\alpha$ have to have the form $\alpha(e^+)=\alpha e^-\tens e^-$ and $\alpha(e^-)=\beta e^+\tens e^+$ in order to commute with functions, for some functions $\alpha,\beta\in A$ (won't refer any more to the map $\alpha$ here). Likewise, to be a bimodule map and lead to a torsion free connection, we must have the form
\begin{gather*}\sigma(e^+\tens e^+)=(1+\gamma)e^+\tens e^+,\quad \sigma(e^+\tens e^-)=a e^+\tens e^-+(1+a)e^-\tens e^+\\
\sigma(e^-\tens e^-)=(1+\delta)e^-\tens e^-,\quad \sigma(e^-\tens e^+)=b e^-\tens e^++(1+b)e^+\tens e^-\end{gather*}
for some functions $\gamma,\delta,a,b$, where we imposed $\wedge(\id+\sigma)=0$. This data corresponds to the connection
\begin{gather*}\nabla e^+=\gamma e^+\tens e^++\alpha e^-\tens e^-+ag,\quad \nabla e^-=\delta e^-\tens e^-+\beta e^+\tens e^++bg\end{gather*}
as the moduli of torsion free bimodule connections. We now impose  the cotorsion equation
\begin{align*} coT_\nabla&=(\extd\tens\id+\id\wedge\nabla)(e^+\tens e^-+e^-\tens e^+)=e^+\nabla e^-+e^-\nabla e^+=\delta e^-+be^++\gamma e^++ae^-=0\end{align*}
which fixes $a=\delta$ and $b=\gamma$. This gives the connection stated as the moduli of bimodule WQLCs. It remains to compute their curvature, 
\begin{align*}R_\nabla(e^+)&=(\extd\tens \id+\id\wedge\nabla)(\alpha e^-\tens e^-+\gamma e^+\tens e^++\delta g)\\
&=\gamma e^+\nabla e^++\alpha e^-\nabla e^-+\del_-\gamma {\rm Vol}\tens e^++\del_+\alpha {\rm Vol}\tens e^-+\extd \delta\wedge g\end{align*}
which computes as stated. Similarly for $R_\nabla(e^-)$. It is also possible to proceed to impose metric compatibility and arrive at Proposition~\ref{proptri} without recourse to the computer result quoted from \cite{MaPac}. \endproof

\subsection{\bf Quantum geometries on a square}\label{sec4pt}  

For $X=\{0,1,2,3\}$, we have 6 connected graphs up to renumbering, each with a unique quantum metric and choices of exterior algebra such as $\Omega_{max},\Omega_{med},\Omega_{min}$. We will not attempt a full classification as we did for $n\le 3$ as this would be a substantial project in its own right. Instead, we focus on the square graph which we take numbered clockwise and which has the merit of being a Cayley graph for two different groups $\Z_4$ and $\Z_2\times \Z_2$, leading to different but natural exterior algebras. They will be different quotients of $\Omega_{min}$. Here  $\Omega^1=\<01,10,12,21,23,32,30,03\>$ in cyclic order, which is 2-dimensional over the algebra. Our general constructions for exterior algebras depend on
\[ \CN_{max}=\CN_{med}=\<012+032,210+230, 123+103, 321+301\>\]
\[ \CN_{min}=\<012+032,210+230, 123+103, 321+301, 010+030,121+101,232+212, 303+323\>\]
while 
\[ \Omega^1\tens_A\Omega^1=\<012,210,123,321,230,032,301,103, 010,030,121,101,232,212, 303,323\>\]
 is spanned by all possible 2-steps. Hence $\Omega^2_{max}=\Omega^2_{med}$ is 12 dimensional as a vector space, in fact free and 3-dimensional over the algebra, while $\Omega_{min}$ is 8-dimensional as a vector space, free and 2-dimensional over the algebra. The unique metric is
\[ g=010+030+121+101+232+212+ 303+323\in \Omega^1\tens_A\Omega^1.\]

\subsubsection{Square with $\Omega(\Z_4)$ calculus} For a natural $\Z_4$-invariant description of the calculus we set
\[ e^+=01+12+23+30,\quad e^-=10+21+32+03\]
for the clockwise and anticlockwise arrows around the square. These are a basis for $\Omega^1$ with relations $e^\pm f=R_\pm(f)e^\pm$ for any $f\in A$, with $(R_\pm f)(i)=f(i\pm 1)$ now referring to $i$ mod 4. In these terms, the inner generator is  $\theta=e^++e^-$ and the quantum metric is $g=e^+\tens e^-+e^-\tens e^+$,  as for the triangle case. This time,  $\Omega_{max}=\Omega_{med}$ is generated over $A$  by $e^\pm$ with the single relation
\[ (e^+){}^2+(e^-){}^2=0,\]
 while $\Omega_{min}$ has the further relation $e^+\wedge e^-+e^-\wedge e^+=0$. In both cases $\extd e^\pm=0$. The standard $\Omega(\Z_4)$ according to the theory of exterior algebras on Cayley graphs\cite{BegMa} has one more relation, i.e. a quotient of $\Omega_{min}$, by setting $(e^\pm)^2=0$ separately, so this is a nontrivial quotient of $\Omega_{min}$. Here 
 \[ \CN_{\Z_4}=\< 012,032, 210,230, 123,103, 321,301, 010+030,121+101,232+212, 303+323 \>.\]
Here, $\Omega^2(\Z_4)$ is again 1-dimensional over the algebra with top `volume' form ${\rm Vol}=e^+\wedge e^-$. It is inner with $\theta=e^++e^-$. We again have a metric inner product of the form $(ij,ji)=\delta_i$ for allowed arrows $j=i\pm 1$ mod 4, and two natural left-invariant lifts given by (\ref{ipm})
as for the triangle. 

\begin{proposition}\label{propsq} For $\Omega(\Z_4)$, the unique quantum metric has  exactly four QLCs  
\begin{align*}\nabla e^+&=\alpha e^-\tens e^-+\alpha\beta(g+e^+\tens e^+),\quad \nabla e^-=\beta e^+\tens e^++\alpha\beta(g+e^-\tens e^-),\quad \alpha,\beta\in\{0,1\)\\
	\sigFun{+}{+} &=(1+\alpha\beta)  \twoForm{+}{+}+ \alpha \twoForm{-}{-} , \quad 
	\sigFun{+}{-} = (1+\alpha\beta)e^-\tens e^++\alpha\beta  \twoForm{+}{-} \\
	\sigFun{-}{-} &=(1+\alpha\beta) \twoForm{-}{-}+\beta \twoForm{+}{+}
,\quad \sigFun{-}{+} = (1+\alpha\beta)e^+\tens e^-+\alpha\beta \twoForm{-}{+}.		\end{align*}
These are all flat.
\end{proposition}
\proof Being a bimodule map requires $\alpha(e^\pm)=0$ in the general construction for a bimodule connection, so these
are classified by $\sigma$ being a bimodule map. If we also impose torsion freeness in the form $\wedge(\id+\sigma)$, we are forced to the form 
\begin{gather*}\sigma(e^+\tens e^+)=\alpha e^-\tens e^-+(1+\gamma) e^+\tens e^+,\quad \sigma(e^-\tens e^-)=(1+\delta) e^-\tens e^-+\beta e^+\tens e^+\\ 
 \sigma(e^+\tens e^-)=ae^-\tens e^++(1+a)e^+\tens e^-,\quad \sigma(e^-\tens e^+)=b e^+\tens e^-+(1+b)e^-\tens e^+\end{gather*}
for some functions $\alpha,\beta,\gamma,\delta, a, b$. This has
\[ \nabla e^+=\alpha e^-\tens e^-+\gamma e^+\tens e^++ (1+a)g,\quad \nabla e^-=\beta e^+\tens e^++\delta e^-\tens e^-+(1+b)g\]
and has too many variables to solve easily for $\nabla g=0$,  but we can first impose the cotorsion equation
\[ coT_\nabla=(\extd\tens \id+\id\wedge\nabla)(e^+\tens e^-+e^-\tens e^+)=e^+\wedge \nabla e^-+e^-\wedge\nabla e^+=0\]
using the relations of the exterior algebra, which forces $b=1+\gamma$ and $a=1+\delta$. We thus have a 4-functional parameter space of bimodule WQLCs
\begin{align}\label{sqwqlc} \nabla e^+&=\alpha e^-\tens e^-+\gamma e^+\tens e^++ \delta g,\quad \nabla e^-=\beta e^+\tens e^++\delta e^-\tens e^-+\gamma g,\\ 
\sigma(e^+\tens e^+)&=\alpha e^-\tens e^-+(1+\gamma) e^+\tens e^+,\quad \sigma(e^-\tens e^-)=(1+\delta) e^-\tens e^-+\beta e^+\tens e^+,\nonumber\\ 
 \sigma(e^+\tens e^-)&=(1+\delta) e^-\tens e^++\delta e^+\tens e^-,\quad \sigma(e^-\tens e^+)=(1+\gamma) e^+\tens e^-+\gamma e^-\tens e^+.\nonumber\end{align}
For QLCs, we now impose metric compatibility
\begin{align*}0&=\nabla(e^+\tens e^-+e^-\tens e^+)\\
&=\alpha e^-\tens e^-\tens e^-+\gamma e^+\tens e^+\tens e^-+\delta g\tens e^- +\beta e^+\tens e^+\tens e^++\delta e^-\tens e^-\tens e^++ \gamma g\tens e^+\\
&\quad +R_+\beta\sigma(e^+\tens e^+)\tens e^++R_+\delta \sigma(e^+\tens e^-)\tens e^-+ R_+\gamma (\sigma(e^+\tens e^+)\tens e^-+\sigma(e^+\tens e^-)\tens e^+)\\
&\quad +R_-\alpha\sigma(e^-\tens e^-)\tens e^-+ R_-\gamma \sigma(e^-\tens e^+)\tens e^++ R_-\delta(\sigma(e^-\tens e^+)\tens e^-+\sigma(e^-\tens e^-)\tens e^+). 
\end{align*}
Substituting the values of $\sigma$ and picking off the coefficients of the tensor products of $e^\pm$ gives 8 equations
\[ \gamma+R_+\gamma + \gamma R_+\gamma= \beta R_-\alpha, \quad \delta+ R_-\delta+ \delta R_-\delta= \alpha R_+\beta,\]
\[ \gamma+R_-\gamma+ \gamma R_-\gamma=\delta R_+\gamma,\quad \delta+R_+\delta+ \delta R_+\delta=\gamma R_-\delta,\]
\[ \delta+ R_-\delta+ \delta R_+\delta=\gamma R_-\delta,\quad \gamma+\gamma R_-\gamma+ R_+\gamma=\delta R_+\gamma,\]
\[ (1+R_-\delta)\beta=(1+\gamma)R_+\beta,\quad(1+R_+\gamma) \alpha= (1+\delta)R_-\alpha.\]
We first deduce that $R_+\gamma=R_-\gamma=\bar\gamma$, say, and $R_+\delta=R_-\delta=\bar\delta$, say. Then $\gamma+\bar\gamma+\gamma\bar\gamma=\delta+\bar\delta+\delta\bar\delta=\gamma\bar\delta=\bar\gamma\delta=a$ for some constant $a$. If $a=0$ then $\gamma=\delta=0$ and if $a=1$ then $\gamma=\delta=1$.  Hence $\gamma=\delta$ is a constant. If $\gamma=\delta=0$ then the last two equations from displayed list tell us that $\alpha,\beta$ are constant and the first two say that $\alpha\beta=0$. If $\gamma=\delta=1$ then $\beta R_-\alpha=1$ from the first equation tells us that $\beta=\alpha=1$ are again constant. Hence we are forced to constant coefficients with $\gamma=\delta=\alpha\beta$ as the only conditions, giving the 4 QLCs stated. That their  curvatures all vanish is a further computation, e.g., 
\[ R_\nabla(e^+)=(\id\wedge\nabla)(\alpha\beta e^+\tens e^++\alpha e^-\tens e^-)=\alpha\beta e^+\wedge\nabla e^++\alpha e^-\wedge\nabla e^-\]
given that the cotorsion already vanishes. We then substitute the value of $\nabla$ to find 0. Similarly for $R_\nabla(e^-)=0$.
\endproof

Therefore, for the square with this calculus, if we want curvature we will have to work more generally with the WQLC's (\ref{sqwqlc}) found during the proof. A short computation for these in the case of constant coefficients gives
\begin{equation}\label{sqcurv}R_\nabla(e^\pm)=(\alpha\beta+\gamma\delta){\rm Vol}\tens e^\pm.\end{equation}
(If the coefficients are not constant then we have derivative terms as in Proposition~\ref{propwqlctri}.) In this case, we have the same form of Ricci tensors ${\rm Ricci}_\pm$ as for the triangle, but now with the factor $\alpha\beta+\gamma\delta$. If we set $\gamma=\delta=0$ then we have the same curvature as in Proposition~\ref{proptri} for the triangle.

\subsubsection{Square with $\Omega(\Z_2\times\Z_2)$ calculus} This has different invariant 1-forms 
\begin{gather*}e^1= (0,0)\to (1,0)+ (1,0)\to (0,0)+(0,1)\to (1,1)+(1,1)\to (0,1),\\
e^2=(0,0)\to (0,1)+ (0,1)\to (0,0)+(1,0)\to (1,1)+ (1,1)\to (1,0), \end{gather*}
where we use Cartesian coordinates for the square, labelled by $\Z_2\times\Z_2$. If we identify this in terms of  our previous labelling of the vertices by $(0,0)=0, (0,1)=1, (1,1)=2, (1,0)=3$ then  using our previous compact notations, this is 
\[e^1=03+30+12+21,\quad e^2=01+10+32+23\]
obeying $e^if=R_i(f)e^i$ for suitable $R_1$ shifting by $(1,0)$ in Cartesian coordinates and $R_2$ by $(0,1)$. The exterior derivative is $\extd f=(\del_i f)e^i$ where $\del_i=R_i+\id$ as we work over $\F_2$. The inner generator and the unique quantum metric on the square in these terms are
\[\theta=e^1+e^2,\quad g=e^1\tens e^1+e^2\tens e^2,\]
while $\Omega_{max}=\Omega_{med}$ is generated over $A$ by $e^i$ with the one relation
\[ e^1\wedge e^2+e^2\wedge e^1=0\]
and $\Omega_{min}$ with the additional relation
\[ (e^1)^2+(e^2)^2=0.\]
The canonical Cayley graph calculus $\Omega(\Z_2\times\Z_2)$ is again a nontrivial quotient of $\Omega_{min}$, with one further relation so that $(e^i)^2=0$ holds separately. This corresponds to 
\[ \CN_{\Z_2\times\Z_2}=\<012+032,210+230, 123+103, 321+301,010, 030,121, 101, 232, 212, 303, 323\>\]
which we see is different from the case of $\Omega(\Z_4)$ previously. Working in $\Omega(\Z_2\times\Z_2)$, there is a unique {\em but not central} volume form and two natural lifts 
\begin{equation}\label{ialt} {\rm Vol}=e^1\wedge e^2=[032+301+ 123+210],\quad i_1({\rm Vol})=e^1\tens e^2,\quad i_2({\rm Vol})=e^2\tens e^1.\end{equation}

\begin{proposition}\label{propsqalt} For $\Omega(\Z_2\times\Z_2)$, the unique quantum metric on the square has exactly 4 QLCs, all flat. Two have constant coefficients,
 \begin{align*} \nabla e^1=\nabla e^2&=\alpha \theta\tens\theta,\quad \alpha\in \{0,1\},\\
  \sigma(e^1\tens e^1)&=\alpha e^2\tens e^2+(1+\alpha) e^1\tens e^1,\quad \sigma(e^2\tens e^2)=(1+\alpha) e^2\tens e^2+\alpha e^1\tens e^1,\\
 \sigma(e^1\tens e^2)&=(1+\alpha) e^2\tens e^1+\alpha e^1\tens e^2,\quad \sigma(e^2\tens e^1)=(1+\alpha) e^1\tens e^2+\alpha e^2\tens e^1\end{align*}
and two have the form
\begin{align*} \nabla e^1&=\theta\tens\theta+\gamma e^1\tens e^1,\quad \nabla e^2=\theta\tens\theta+(1+\gamma)e^2\tens e^2\\
\sigma(e^1\tens e^1)&=\gamma e^1\tens e^1+ e^2\tens e^2,\quad \sigma(e^1\tens e^2)=e^1\tens e^2\\
\sigma(e^2\tens e^2)&=e^1\tens e^1+ (1+\gamma)e^2\tens e^2,\quad \sigma(e^2\tens e^1)=e^2\tens e^1,
\end{align*}
where $\gamma$ is a function that alternates between 0 and 1 as we go around the square (hence determined by its value at one vertex). 
\end{proposition}
\proof Being a bimodule map forces $\alpha(e^i)=0$ since $i=1,2$ correspond to $(0,1),(1,0)$ in the group. Similarly, to be a bimodule map and give a torsion free connection, we are forced to
\begin{align*}\sigma(e^1\tens e^1)=\alpha e^2\tens e^2+\gamma e^1\tens e^1,\quad \sigma(e^2\tens e^2)=\delta e^2\tens e^2+\beta e^1\tens e^1,\\
 \sigma(e^1\tens e^2)=ae^2\tens e^1+(1+a)e^1\tens e^2,\quad \sigma(e^2\tens e^1)=b e^1\tens e^2+(1+b)e^2\tens e^1\end{align*}
for functions $\alpha,\beta,\gamma,\delta, a, b$. This coincidentally has the same form as at the start of the proof of Proposition~\ref{propsq} and the connection therefore has the same form, 
\begin{align*} \nabla e^1=\alpha e^2\tens e^2+(1+\gamma)e^1\tens e^1+ (1+a)(e^1\tens e^2+e^2\tens e^1),\\
 \nabla e^2=\beta e^1\tens e^1+(1+\delta)e^2\tens e^2+(1+b)(e^1\tens e^2+e^2\tens e^1).\end{align*}
The metric, however, now has a different form, so the cotorsion is not parallel. This time
\[ coT_\nabla=(\extd\tens\id+ \id\wedge\nabla(e^1\tens e^1+e^2\tens e^2)=e^1\nabla e^1+e^2\nabla e^2=0\]
 now forces $b=1+\alpha$ and $a=1+\beta$. We again have a 4-functional parameter space of  bimodule WQLCs, this time
 \begin{align}\label{sqwqlcalt} \nabla e^1&=\alpha e^2\tens e^2+(1+\gamma)e^1\tens e^1+ \beta (e^1\tens e^2+e^2\tens e^1),\nonumber\\
 \nabla e^2&=\beta e^1\tens e^1+(1+\delta)e^2\tens e^2+\alpha (e^1\tens e^2+e^2\tens e^1),\\
 \sigma(e^1\tens e^1)&=\alpha e^2\tens e^2+\gamma e^1\tens e^1,\quad \sigma(e^2\tens e^2)=\delta e^2\tens e^2+\beta e^1\tens e^1,\nonumber\\
 \sigma(e^1\tens e^2)&=(1+\beta) e^2\tens e^1+\beta e^1\tens e^2,\quad \sigma(e^2\tens e^1)=(1+\alpha) e^1\tens e^2+\alpha e^2\tens e^1.\nonumber\end{align}
 For QLCs, we now impose metric compatibility
 \begin{align*}0&=\nabla(e^1\tens e^1+e^2\tens e^2)\\
 &=\alpha e^2\tens e^2\tens e^1+(1+\gamma)e^1\tens e^1\tens e^1+\beta(e^1\tens e^2+e^2\tens e^1)\tens e^1\\
 &\quad +\beta e^1\tens e^1\tens e^2+(1+\delta)e^2\tens e^2\tens e^2+\alpha(e^1\tens e^2+e^2\tens e^1)\tens e^2\\
 &\quad + (R_1\beta \sigma(e^1\tens e^2)+R_2\beta\sigma(e^2\tens e^1)+R_2\alpha\sigma(e^2\tens e^2)+(1+R_1\gamma)\sigma(e^1\tens e^1))\tens e^1\\
 &\quad +(R_2\alpha\sigma(e^2\tens e^1)+R_1\alpha\sigma(e^1\tens e^2)+R_1\beta\sigma(e^1\tens e^1)+(1+R_2\delta)\sigma(e^2\tens e^2))\tens e^2\end{align*}
 giving for the coefficients of the tensor powers of the $e^i$ the 8 equations
 \[ \gamma R_1\gamma=1+\beta R_2\alpha,\quad \delta R_2\delta=1+\alpha R_1\beta\]
 \[ \delta R_2\alpha=\alpha R_1\gamma,\quad \gamma R_1\beta=\beta R_2\delta\]
 \[ \beta+R_2\beta+ \beta R_1\beta =\alpha R_2\beta,\quad \alpha+R_1\alpha+\alpha R_2\alpha=\beta R_1\alpha\]
 \[ \beta+R_1\beta+\beta R_1\beta=\alpha R_2\beta,\quad \alpha+R_1\alpha+\alpha R_1\alpha=\beta R_1\alpha\]
 From the last two lines, we conclude that $R_1\alpha=R_2\alpha=\bar\alpha$, say, and $R_1\beta=R_2\beta=\bar\beta$. We then deduce that $\alpha+\bar\alpha+\alpha\bar\alpha=\beta+\bar\beta+\beta\bar\beta=\alpha\bar\beta=\bar\alpha\beta=a$, a constant. If $a=0$ then we are forced to $\alpha=\beta=0$ and if $a=1$ then we are forced to $\alpha=\beta=1$. We conclude that $\alpha=\beta$ is a constant. If $\alpha=\beta=0$ then the first of our equations tells us that $\gamma=\delta=1$. This is one of our four solutions. If $\alpha=\beta=1$ then we are forced to $\delta=R_1\gamma$ and are only left to solve  $R_2R_1\gamma=\gamma$ and $\gamma R_1\gamma=0$ for $\gamma$. This is solved by $\gamma(0,0)=\gamma(1,1)=a$ and $\gamma(1,0)=\gamma(0,1)=b$ such that $ab=1$, which is a further three solutions. All four are contained in the definitions $\alpha=1+ab$, $\delta=R_1\gamma$ and $\gamma$ as specified. Here $a=b$ gives two solutions with constant coefficients $\gamma=\delta=1+\alpha$, while $b=1+a$ gives two solutions with $\alpha=1$, nonconstant coefficients $\gamma$ as stated, and $\delta=1+\gamma$.  It is easy to see that the constant solutions are flat (they have a similar structure to the $\alpha=\beta$ case of Proposition~\ref{propsq}). For the nonconstant solutions we have
 \begin{align*} R_\nabla(e^1)&=(\extd\tens\id+\id\wedge\nabla)(\theta\tens\theta+\gamma e^1\tens e^1)=\extd\gamma\wedge e^1\tens e^1+\theta \wedge\nabla\theta+\gamma e^1\wedge \nabla e^1\\
 &=\gamma {\rm Vol}\tens e^1+{\rm Vol}\tens e^2+\gamma{\rm Vol}\tens e^2+R_1\gamma {\rm Vol}\tens e^2+\gamma {\rm Vol}\tens\theta=0\end{align*}
using that the calculus is inner and that $\del_1\gamma=1$ for our particular $\gamma$. Similarly for $R_\nabla(e^2)=0$.  \endproof

If we want curvature, we can look more generally to the WQLCs (\ref{sqwqlcalt}) found during the above proof. A short computation for these in the case of constant coefficients gives
\begin{equation}\label{sqcurvalt} R_\nabla(e^1)=(\alpha\gamma+\beta\delta){\rm Vol}\tens e^2,\quad R_\nabla(e^2)=(\alpha\gamma+\beta\delta){\rm Vol}\tens e^1.\end{equation}
(If the coefficients are not constant then we have derivative terms as usual.) These WQLC curvatures have two natural Ricci tensors
\[ {\rm Ricci}_1=(\alpha\gamma+\beta\delta)e^2\tens e^2,\quad {\rm Ricci}_2=(\alpha\gamma+\beta\delta)e^1\tens e^1,\quad S_1=S_2=\alpha\gamma+\beta\delta\]
according to the two lifts (\ref{ialt}), with sum ${}^2{\rm Ricci}=(\alpha\gamma+\beta\delta)g$.

\subsection{Quantum geometry on an $n$-gon with $n\ge 5$}\label{secpoly}

For larger $n$, the number of possible connected graphs explode rapidly and there are many possibilities even if we focus on Cayley graphs, namely as the product of cyclic groups according to the prime factorisation of $n$. We have seen this already for $n=4$ with its two factorisations. Here we focus just on the $n$-gon case, as the Cayley graph for $\Z_n$ with its standard generator. We have already  covered $\Z_3$  in Propositions~\ref{proptri} and $\Z_4$ in Proposition~\ref{propsq}. For general $n\ge 5$, we similarly number the vertices in sequence clockwise and define
\[ e^+=\sum_i i\to i+1,\quad  e^-=\sum_i i\to i-1\]
with bimodule relations $e^\pm f=(R_\pm f)e^\pm$, where $(R_\pm f)(i)=f(i\pm 1)$ for $i$ mod $n$. We have
\[ \CN_{max}=\CN_{med}=\<i\to i+1\to i+2, i\to i-1\to i-2\>,\]
\[\CN_{min}=\<i\to i+1\to i+2, i\to i-1\to i-2, (i\to i+1\to i)+( i\to i-1\to i)\> \]
leading similarly to $\Omega_{max}=\Omega_{med}$ and $\Omega_{min}$ free with 2-forms $2$-dimensional and $1$-dimensional respectively over $A=\F_2(\Z_n)$. In terms of the $e^\pm$, the unique quantum metric over $\F_2$ is 
\[ g=\sum_i (i\to i+1\to i)+(i\to i-1\to i)=e^+\tens e^-+e^-\tens e^+\]
while $\Omega_{max}=\Omega_{med}$ has the relations $(e^\pm)^2=0$ but no other relations among them, and $\extd e^\pm=0$. By contrast, $\Omega_{min}=\Omega(\Z_n)$, the canonical Cayley graph exterior algebra, with the additional relation $e^+\wedge e^-+e^-\wedge e^+=0$ . In this case, we have
\[  {\rm Vol}=e^+\wedge e^-=[\sum_i i\to i+1\to i],\quad i_\pm({\rm Vol})=e^\pm\tens e^\mp\]
as the natural top or volume form (the unique invariant one) and two natural lifts (the only invariant ones over $\F_2$). This description of $\Omega_{med}$ and $\Omega_{min}$ also applied to $n=3$. 

\begin{proposition} Let $n\ge 5$. For  $\Omega_{min}=\Omega(\Z_n)$, the unique quantum metric has a unique QLC, $\nabla e^\pm=0$. This has $\sigma={\rm flip}$ on the $e^\pm$ and is flat.
\end{proposition}
\proof As for $n=4$, to be bimodule map we need $\alpha(e^\pm)=0$. The form of a bimodule map $\sigma$ after imposing zero torsion by $\wedge(\id+\sigma)=0$ is now
\[ \sigma(e^+\tens e^+)=(1+\alpha) e^+\tens e^+,\quad \sigma(e^+\tens e^-)=(1+\delta) e^-\tens  +\delta e^-\tens e^+\]
\[ \sigma(e^-\tens e^-)=(1+\beta)e^-\tens e^-,\quad \sigma(e^-\tens e^+)=(1+\gamma)e^+\tens e^-+\gamma e^-\tens e^+\]
so that there is a 4-functional parameter space of torsion free bimodule connections  
\[ \nabla e^+=\alpha e^+\tens e^++ \delta g,\quad \nabla e^-=\beta e^-\tens e^-+ \gamma g.\]
We now impose $coT_\nabla=e^+\wedge \nabla e^-+e^-\wedge\nabla e^+=0$ as in the proof of Proposition~\ref{propsq}, which forces $\delta=\beta$ and $\gamma=\alpha$. Thus we have a 2-functional parameter space of bimodule WQLCs
\begin{gather} \nabla e^+=\alpha e^+\tens e^++ \beta g,\quad \nabla e^-=\beta e^-\tens e^-+ \alpha g,\nonumber\\
\label{polywqlc} \sigma(e^+\tens e^+)=(1+\alpha) e^+\tens e^+,\quad \sigma(e^+\tens e^-)=(1+\beta) e^-\tens e^+  +\beta e^-\tens e^+,\\
 \sigma(e^-\tens e^-)=(1+\beta)e^-\tens e^-,\quad \sigma(e^-\tens e^+)=(1+\alpha)e^+\tens e^-+\alpha e^-\tens e^+. \nonumber\end{gather}
For QLCs, we now impose metric compatibility
\begin{align*}0&=\nabla(e^+\tens e^-+e^-\tens e^+)\\
&=\alpha e^+\tens e^+\tens e^-+\beta g\tens e^-+ \beta e^-\tens e^-\tens e^++\alpha g\tens e^+\\
&\quad + R_+\alpha(\sigma(e^+\tens e^+)\tens e^-+ \sigma(e^+\tens e^-)\tens e^+)+ R_+\beta \sigma (e^+\tens e^-)\tens e^-\\
&\quad + R_-\beta (\sigma(e^-\tens e^+)\tens e^-+ \sigma(e^-\tens e^-)\tens e^+) + R_-\alpha \sigma (e^-\tens e^+)\tens e^+.
\end{align*}
Substituting the values of $\sigma$  again gives 8 equations. The coefficients of $e^+\tens e^+\tens e^-$ and $e^-\tens e^-\tens e^+$ tell us that 
\[ \alpha+R_+\alpha+\alpha R_+\alpha=0,\quad \beta+ R_-\beta+ \beta R_-\beta=0\]
which imply that $\alpha=\beta=0$. Thus, the only QLC is $\nabla e^\pm=0$, with $R_\nabla=0$. \endproof

As with $n=4$, if we want curvature then we have to turn to WQLCs. These were found during the proof in (\ref{polywqlc}) and in the case of constant coefficients, have
\begin{equation}\label{polycurv} R_\nabla e^+=\id\wedge\nabla(\alpha e^+\tens e^+)=\alpha e^+\wedge\nabla e^+=\alpha\beta{\rm Vol}\tens e^+\end{equation}
and similarly $R_\nabla e^-=\alpha\beta{\rm Vol}\tens e^-$. This is the same as for the triangle  in Proposition~\ref{proptri} and the two Ricci tensors for the lifts $i_\pm$ are also the same, by the same calculation. If we allowed non-constant coefficients for the WQLCs then we would have derivative terms in the curvature as in Proposition~\ref{propwqlctri}.
 
\subsection{Polygon geometry in subset Boolean form}\label{secsub}

We have already exhibited the $0-1$ graph and the $0-1-2$ graph geometries in subset notation. For the polygon geometry, we made extensive use of the left-invariant basis method and here our results require more involved translation to the subset description.  We cover only $n\ne 4$ so that we are always working with $\Omega_{min}=\Omega(\Z_n)$ in agreement with the Cayley graph calculus for this group. 

\subsubsection{Subset version of triangle}\label{subsettri}

 Here $\Omega^1=P(\{01,10,12,21,20,02\})$ containing $\theta=\{01,10,12,21,20,02\}$, 
 \[ \Omega^1\tens_A\Omega^1=P(\{020,010,012,021,101,121,120, 102, 212, 202,201,210\})\]
 containing
 \[ g=010+020+101+121+202+212\]
 and the values of the QLCs. From the Leibniz rule on products $\delta_ie^\pm$ and the given  $\nabla e^\pm$, one can deduce from Proposition~\ref{proptri} for $\Omega_{min}=\Omega(\Z_3)$:

(i) For $\alpha=\beta=0$ 
\begin{gather*}\nabla\{ 01\}=\{020, 201, 101, 012 \},\quad \nabla\{12\}=\{101,012,212,120\},\quad \nabla\{20\}=\{212,120,020,201\},\\
 \nabla\{10\}=\{121,210,010,102\},\quad\nabla\{21\}=\{202,021,121,210\},\quad  \nabla\{02\}=\{010,102, 202, 021\},\end{gather*}
Others are given by $\oplus$, for example $\oplus$ of each row gives $\nabla e^\pm=\emptyset$. 

(ii) For $\alpha=1,\beta=0$ we $\oplus$ an extra term to (i) for increasing arrows only, according to the bimodule map (also denoted $\alpha$):
\[ \alpha(\{ 01\})=\{021\},\quad \alpha(\{12\})=\{102\},\quad \alpha(\{20\})=\{210\}.\]

(iii) For $\alpha=0,\beta=1$ we $\oplus$ an extra term to (i) for decreasing arrows only, according to the bimodule map
\[  \alpha(\{10\})=\{120\},\quad \alpha(\{21\})=\{ 201\},\quad \alpha(\{02\})=\{012\}.\]

(iv) $\alpha=1,\beta=1$ we $\oplus$ both  (ii) and (iii) as applicable to the connection in (i). This case  has curvature given by $R_\nabla={\rm Vol}\tens$, so for example
\[ \nabla\{01,20\}=\{101, 012, 021, 212, 120, 210\},\quad R_\nabla\{01,20\}=\{0101, 2020\}.\]

\subsubsection{Subset version of $n$-gon}

Here, we translate in general how the unique $\nabla e^\pm=0$ QLC on the polygon with $\Omega_{min}=\Omega(\Z_n)$ looks in terms of subsets for $n\ge 5$. These results also apply to $n=3$ for this QLC, just this is no longer unique. In the general case, we do not have the luxury of the compact notation and write the arrows explicitly. Thus
\[ e^+ =\{i\to i+1\ |\ i=0,\cdots,n-1\}; \quad e^- =\{i\to i-1\ |\ i=0,\cdots,n-1\}\]
as elements of $\Omega^1=P({\rm Arr})$ where ${\rm Arr}=\{i\to i\pm 1\}$ is partitioned into singleton sets
\[ {}_i{\rm Arr}_{i+1}=\{i\to i+1\},\quad {}_i{\rm Arr}_{i-1}=\{i\to i-1\}\]
going respectively clockwise and anticlockwise around the $n$-gon numbered clockwise by $\Z_n$. This $e^+$ is the union of the first and $e^-$ of the second. Moreover, a general subset of arrows can be expressed in the form  $\omega=(a_+\cap e^+)\cup (a_-\cap e^-)$ for some $a_\pm\subseteq P(X)$ which can be recovered from $\omega$ by
\begin{equation}\label{apm} a_\pm=\{{\rm tails\ of\ } \omega\cap e^\pm\}\subseteq X.\end{equation}

Next, $\Omega^1\tens_A\Omega^1=P({\rm Arr}^{(2)})$ with the $4n$ 2-step arrows ${\rm Arr}^{(2)}$ partitioned into subsets
\[ {}_i{\rm Arr}^{(2)}_i=\{i\to i+1\to i, i\to i-1\to i\},\quad  {}_i{\rm Arr}^{(2)}_{i+2}=\{i\to i+1\to i+2\},\quad {}_i{\rm Arr}^{(2)}_{i-2}=\{i\to i-1\to i-2\}.\]
The canonical $\CN_{min}$  sets all of these subsets of arrows  as well as all their unions to zero (in the sense of an equivalence relation on $P({\rm Arr}^{(2)})$). Then $\Omega^2$ is $n$-dimensional over $\F_2$ with every element represented as $a\cap {\rm Vol}$ where
\begin{equation}\label{volpoly}  {\rm Vol}=\cup_i\{i\to i+ 1\to i\}\end{equation}
in the quotient and  in agreement with the Cayley graph construction. The metric as a subset is 
\begin{equation}\label{metpoly} g=\cup_i \left({}_i{\rm Arr}^{(2)}_i\right)=\cup_i\{i\to i\pm 1\to i\}\end{equation}
and we see that $\wedge(g)=0$ in the quotient as the two entries for each $i$ are equivalent with respect to $\CN_{min}$.   Later on, for the Ricci tensor, we will need a lift $\Omega^2\to \Omega^1\tens_{P(X)}\Omega^1$ and we have two natural $\Z_n$-invariant ones
\begin{equation}\label{ipmpoly} i_+({\rm Vol})= \cup_i\{i\to i+1\to i\},\quad i_-({\rm Vol})=\cup_i\{i\to i-1\to i\}\end{equation}
amounting to two halves of the metric.  

We next see how the trivial connection $\nabla e^\pm=0$ given by  $\sigma=${\rm flip} on the generators and the bimodule map $\alpha=0$, is verified to be a QLC in terms of $P(X)$. Here,  $\sigma: P({\rm Arr}^{(2)})\to P({\rm Arr}^{(2)})$ is given elementwise on subsets by the maps 
\begin{equation}\label{sigflip} \sigma|_{{}_i{\rm Arr}^{(2)}_i}={\rm swap},\quad  \sigma|_{{}_i{\rm Arr}^{(2)}_{i+2}}=\id,\quad  \sigma|_{{}_i{\rm Arr}^{(2)}_{i-2}}=\id\end{equation}
where {\rm swap} gives the other element of the 2-element set.  Clearly $\wedge(\id+\sigma)=0$ since in the quotient the swap in the first map has no effect. So $\nabla$ defined by $\sigma$ and the map $\alpha=0$ in (\ref{nablatheta}) is torsion free. For metric compatibility, we have
\[ g\tens\theta=\cup_p \{p\to p\pm 1\to p\to p\pm'1\},\quad \theta\tens g=\cup_p \{p\pm1\to p\to p\pm'1\to p\},\]
where $\pm'$ are independent so that each set has 4 elements. Applying $\sigma$ as in (\ref{sigflip}),  one readily sees that as sets $\sigma_{12}\sigma_{23}(g\tens\theta)=\theta\tens g$ so that $\nabla$  is metric compatible by (\ref{metcompsig}) and hence a QLC.

\begin{figure}
\[ \includegraphics[scale=0.5]{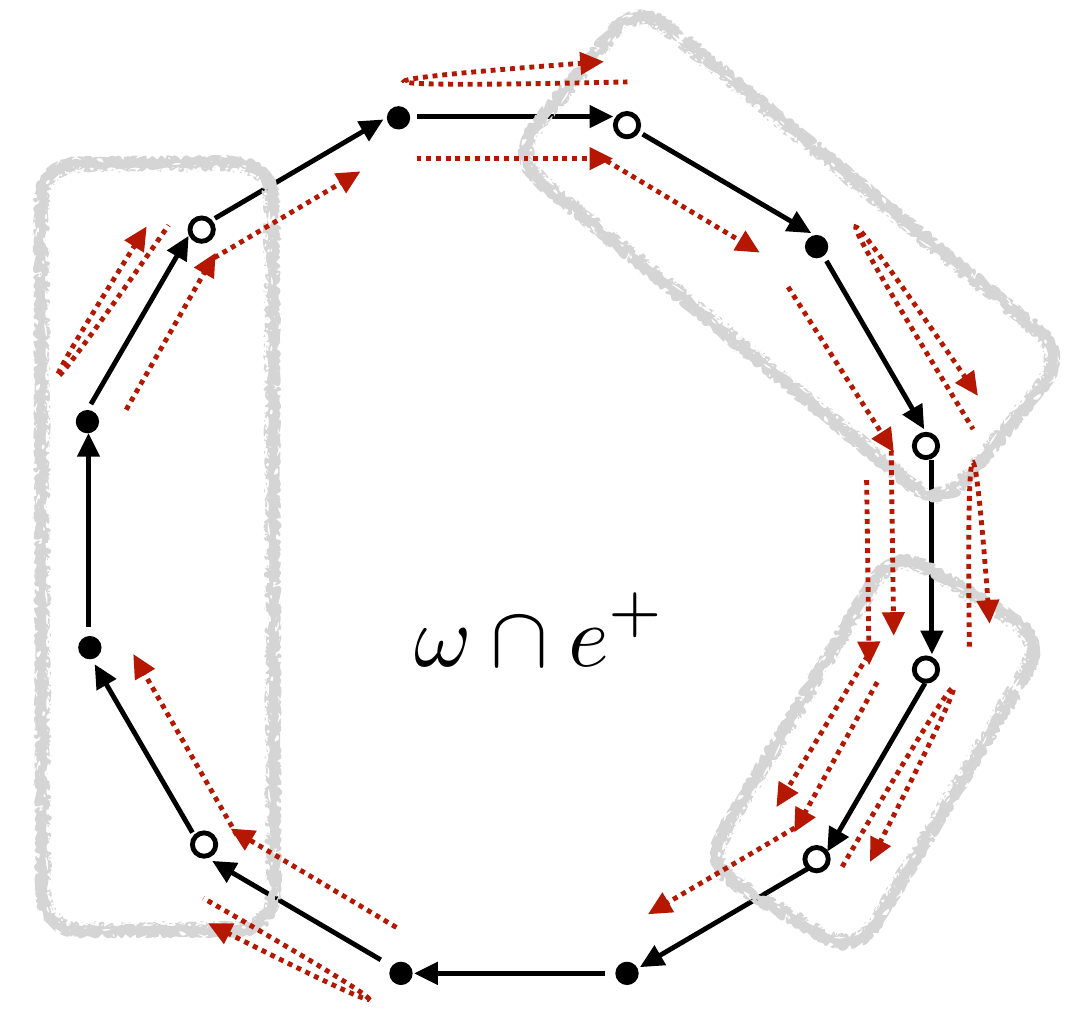}\]
\caption{\label{fignabla} Half of the trivial connection $\nabla\omega$. Here the boxes are $\omega\cap e^+$  and the open circles are the boundary points $\del_+\omega$. Each of these contribute two 2-steps $e^+_i$ to $\nabla\omega$ as shown dashed. The other half of $\nabla$ is the same construction applied to $\omega\cap e^-$ with boundary points $\del_-\omega$ contributing $e^-_i$. }
\end{figure}

To see what $\nabla$ looks like, we compute for example
\begin{align*} \nabla \{i\to i+1\}&=\{i\pm 1\to i\to i+1\}\oplus\sigma(\{i\to i+1\to i+2,i\to i+1\to i\})\\
&=\{i\to i-1\to i, i-1\to i\to i+1, i+1\to i\to i+1, i\to i+1\to i+2\}\\
&=e^+_i\cup e^+_{i+1},\end{align*}
where we define  
\[ e^+_i= \{i\to i-1\to i,i-1\to i\to i+1\},\quad e^-_i=\{ i\to i+1\to i,i+1\to i\to i-1\}\]
and use $e^-_i$ for a similar result $\nabla\{i-1\to i\}=e^-_{i-1}\cup e^-_i$. For the general case, let $\omega\subseteq {\rm Arr}$ and 
\[ \del_+\omega=\{i\in {\rm tails\ of\ }\omega_+ \ | i-1\notin {\rm tails\ of\ }\omega_+\}\cup \{i\in {\rm heads\ of\ }\omega \ | i+1\notin {\rm heads\ of\ }\omega_+\},\]
\[ \del_-\omega=\{i\in {\rm tails\ of\ }\omega_- \ | i+1\notin {\rm tails\ of\ }\omega_-\}\cup \{i\in {\rm heads\ of\ }\omega \ | i-1\notin {\rm heads\ of\ }\omega_-\},\]
where $\omega_\pm=\omega\cap e^\pm$ are the increasing/decreasing arrows of $\omega$. Then 
\begin{equation}\label{trivcon}\nabla: P({\rm Arr})\to P({\rm Arr}^{(2)}),\quad  \nabla\omega=(\cup_{i\in\del_+\omega} e^+_i)\cup (\cup_{i\in \del_-\omega}e^-_i)\end{equation}
as depicted in Figure~\ref{fignabla}.  All components here are disjoint. One can check that this general description reduces for $n=3$ to the explicit formulae in case (i) of Section~\ref{subsettri}. 

\subsection{De Morgan dual connection on the polygon.}\label{secdem}

For the polygon with $n\ge 4$, the trivial connection $\nabla$ of the preceding section has $\nabla(\omega\oplus\bar\omega)=\nabla\theta=\emptyset$ for any $\omega\subseteq{\rm Arr}$, hence $\nabla\bar\omega=\nabla\omega$ and the de Morgan dual connection
\[ \bar\nabla\omega:=\overline{\nabla\bar\omega}=\nabla\omega\oplus{\rm Arr}^{(2)}.\]
One also has $\nabla\bar\omega=\nabla\omega$ and hence the same conclusion $\bar\nabla\omega=\nabla\omega\oplus{\rm Arr}^{(2)}$ for the trivial $\alpha=\beta=0$ QLC for $n=4$ in Proposition~\ref{propsq}. 
In what follows, we now focus on the curved QLC found on the $n$-gon for $n=3$, the $\alpha=\beta=1$ case of Proposition~\ref{proptri}. 

\begin{lemma}\label{lemdemqlctri} For the triangle, the curved $\alpha=\beta=1$ QLC has de Morgan dual connection
\[ \bar\nabla \omega=\nabla\omega\oplus g\]
for all $\omega\subset{\rm Arr}$, where $g$ is the quantum metric. 
\end{lemma}
\proof  If we denote the trivial connection by $\nabla_0e^\pm=0$ then the $\alpha=\beta=1$ connection can be written as
\[ \nabla\omega=\nabla_0\omega\cup\{i\to i-1\to i-2\ |\ i\in a_+\}\cup \{i\to i+1\to i+2\ |\ i\in a_-\},\]
where $a_\pm$ are defined from $\omega$ by (\ref{apm}). Denoting this dependence explicitly, we observe that $a_\pm(\bar\omega)=\overline{a_\pm(\omega)}$, which combined with our observation that $\nabla_0\bar\omega=\nabla_0\omega$, leads to
\[ \nabla\bar\omega=\nabla\omega\oplus  \{i\to i+1\to i-1, i\to i-1\to i+1\ |\ i=0,1,2\}=\nabla\omega\oplus\CN_{med},\]
which in turn tells us that the corresponding connection on $\bar P(X)$ is
\[ \bar\nabla\omega:=\overline{\nabla\bar\omega}=\nabla\bar\omega\oplus{\rm Arr}^{(2)}= \nabla\omega\oplus g,\]
but viewed in $\bar P({\rm Arr}^{(2)})$. \endproof

In terms of left-invariant 1-forms,  this is $\bar\nabla e^\pm=\nabla e^\pm+ g=e^\pm\tens e^\pm + \theta\tens\theta$. 
For the curvature, we first recall that the $\Omega_{min}$ calculus is given by quotienting out
\[ \CN_{min}=\{\{010,020\}, \{121,101\}, \{202,212\}, \{012\},\{120\},\{201\}, \{210\},\{021\},\{102\}\}\]
in the subset form. This defines an equivalence relation on $P({\rm Arr}^{(2)})$ where $\omega\sim\eta$ for  $\omega,\eta\subseteq {\rm Arr}^{(2)}$ if $\omega\oplus\eta$ is a union of some of the subsets in the collection $\CN_{min}$. In particular,
\[ {\rm Vol}\sim \{010,121,202\}\oplus{\rm Arr}^{(2)}.\]
implies from the second form that 
\[\overline{\rm Vol}=\{010,121,202\}={\rm Vol},\]
but now taken modulo $\bar \CN_{min}$. The latter is defined in the same way as $\CN_{min}$ but with $\bar\oplus$ of the  ${}_p{\rm Arr}^{(2)}_q$ subsets. We can check that $\overline{\rm Vol}$ here is well-defined. For example, if we took ${\rm Vol}\sim \{010,121,202,120\}$ then using the second form,
\[ \overline{\rm Vol}=\overline{\{010,121,202\}\oplus\{120\}}=\{010,121,202\}\bar\oplus \{120\}\]
which is equivalent with respect to  $\bar \CN_{min}$. 

Since $\bar{\ }$ is a diffeomorphism, the connection $\bar\nabla$ must have the dual form to the curvature of $\nabla$, so 
\[ R_{\bar\nabla}(\omega)=\overline{{\rm Vol}\tens_{P(X)}\bar\omega}= \overline{\rm Vol}\tens_{\bar P(X)}\omega\]
with the left factor now understood modulo $\bar \CN_{min}$. As a check on our entire dual formalism, we now verify this curvature directly on $\omega=\{01\}$ as follows. During calculations, $\tens$ means $\tens_{P(X)}$ while $\bar\tens=\tens_{\bar P(X)}$. We start with
\begin{align*} \bar\nabla\{01\}&=\{020,201\}\oplus g=\{010,121,012,202,212,201\}=\overline{\{020,101,120,102,210\}}\\
&=\overline{\{02,10,12,21\}\tens \{20,01,02,10\}}=\overline{\{02,10,12,21\}}\bar\tens \overline{\{20,01,02,10\}}\end{align*}
\begin{align*}&\kern-10pt (\bar\extd \overline{\{02,10,12,21\}})\bar\tens\overline{\{20,01,02,10\}}=\{201,021,121,212,210,102\}\bar\tens\overline{\{20,01,02,10\}}\\
&=\overline{\{ 020,010,202,101,120,012\}}\bar\tens\overline{\{20,01,02,10\}}\\
&=\overline{\{ 020,010,202,101,120,012\}\tens\{20,01,02,10\}}\\ &=\overline{\{0202,0201,2020,0101,0102,1010,1202,1201,0120\}}\end{align*}
where $\bar\extd$ on a subset of arrows is all 2-steps wholly in our out of the given subset of arrows. This has a similar form to $\bar\extd$ on subsets of vertices but for an induced graph on ${\rm Arr}$ whereby two arrows have a higher-level arrow if they concatenate. For example, $20\Rightarrow 01$ (since they form a 2-step) contributes the 2-step 201. Its result should be understood modulo $\bar\CN_{min}$. On the other factor of the output of $\bar\nabla$ we compute
\begin{align*}&\kern-10pt\overline{\{02,10,12,21\}}\bar\tens\bar\nabla \overline{\{20,01,02,10\}}=\overline{\{02,10,12,21\}}\bar\tens\overline{\nabla\{20,01,02,10\}}\\
&=\overline{\{02,10,12,21\}}\bar\tens\overline{\{101,121,202,212\}}=\overline{\{02,10,12,21\}\tens \{101,121,202,212\}}\\
&=\overline{\{0202,0212,1202,1212,2101,2121\}} \end{align*}
again with the first 2 steps modulo $\bar\CN_{min}$ (which we denote $\bar\wedge$). The curvature computed in $\bar P(X)$ is the $\bar\oplus$ of these two results:
\begin{align*}&R_{\bar\nabla}\{01\}=(\bar\extd\bar\tens\id \bar\oplus \id\bar\wedge\bar\nabla)\bar\nabla\{01\}\\
&=\overline{\{ 0202,0201,2020,0101,0102,1010,1202,1201,0120  \}\oplus\{  0202,0212,1202,1212,2101,2121 \}}\\
&=\overline{\{ 0101,0102, 0120 , 0201, 0212,1010,1201,1212, 2020, 2101,2121  \}}\\
&=\overline{\{ 0102,1210,1212, 2020, 2121  \}}=\overline{\{010,121,202\}\tens\{10,02,20,12,21\}}=\overline{\rm Vol}\bar\tens\{01\}, 
\end{align*}
where in the 4th equality we use the relations of $\CN_{min}$ inside the overline to simplify.

\section{Generalised de Morgan duality over $\F_2$}\label{secbarA} 

While de Morgan duality is natural in a Boolean context, here we extend it to any unital algebra $A$ over $\F_2$ based on our point of view of Boolean algebras as a special case. This would be rather unusual coming at it from the point of view of noncommutative differential geometry but helps to explore its geometric significance. We define `complementation' as  the  bijection 
\[ A\to A,\quad a\mapsto \bar a=1+a,\]
which we view as an isomorphism of $A$ with a new algebra structure on $A$, denoted $\bar A$, with new product and addition 
\[ a\bar\cdot b= ab+a+b,\quad a\bar+ b=a+b+1.\]

\begin{lemma} $\bar A$ with the above product and addition is again a unital algebra over $\F_2$ with $\bar 1=0$ and $\bar 0=1$. Moreover, it obeys $a\bar\cdot a=a^2$ so the new algebra is Boolean if and only if the initial one is. \end{lemma}
\proof This involves checking all the axioms of an algebra. For example, 
\begin{align*} (a\bar\cdot b)\bar+(a\bar\cdot c)&= 1+a\bar\cdot b+ a\bar\cdot c=1+a+b+ab+a+c+ac\\ &=a+(1+b+c)+a+ab+ac= a\bar\cdot (1+b+c)=a\bar\cdot(b\bar+c)\\
a\bar\cdot(b\bar\cdot c)&=a+(b\bar\cdot c)+a(b\bar\cdot c)=a+b+c+bc+ab+ac+abc=(a\bar\cdot b)\bar\cdot c,\end{align*}
where the last step is similar to the preceding ones but in reverse. We also have $1\bar+ a=1+1+a=a$ over $\F_2$ while $0\bar\cdot a=0+a+0a=a$, which agrees with $1=\bar 0$ and $0=\bar 1$. As a check, we then have $1\bar\cdot a=1+a+a=1$ over $\F_2$, which is $\bar 0\bar\cdot a=\bar 0$. \endproof

This restricts correctly to the atomic Boolean case of $A=\F_2(X)$. From the point of view of characteristic functions, $\bar{\chi_a}=\chi_{\bar a}$ for $a\subseteq X$ and
\[ \chi_a\bar\cdot\chi_b=\chi_{a\cup b},\quad \chi_a\bar+\chi_b=\chi_{\overline{a\oplus b}}=\chi_{a\bar\oplus b}.\] 
Or working directly on $P(X)$, 
\[ a\bar\cdot b=a\oplus b\oplus (a\cap b)=((a\cup b)\cap\overline{a\cap b})\cup (a\cap b)= a\cup b\cup (a\cap b)=a\cup b\]
 \[a\bar+ b= X\oplus (a\oplus b)=\overline{a\oplus b}=a\bar\oplus b\] 
recovers the expected algebra structure of $\bar P(X)$. 
 
In our more algebraic language, however, $A=\F_2[\delta_i]/\< \delta_i\delta_j-\delta_{ij}\delta_j\>$ for $i\in X$, while in $\bar A$, their product obeys
 \[ \delta_i\bar\cdot\delta_j\bar+\delta_i\bar+\delta_j\bar+\bar 1\bar\plus\delta_{ij}(\bar 1\bar\plus\delta_i)=\delta_i+\delta_j+\delta_{ij}\delta_i+1+\delta_i+1+\delta_j+1+0+\delta_{ij}(1+1+0+\delta_i)=1=\bar 0.\]
We used here that if $\mu=0,1$ then $\bar+\mu a$ in $\bar A$ means $\bar+\bar 0 =+0$ in $A$ if $\mu=0$ and $+1+a$ in $A$ if $\mu=1$, i.e.  $\bar+\mu a=+\mu(1+a)$ in $A$. We obtained just the relations in $A$ of the complementary projectors $\eps_i:=1+\delta_i$ obeying $\eps_i\eps_j=1+\eps_i+\eps_j+\delta_{ij}(1+\eps_i)=0$. Thus one can also view de Morgan duality as a change of variables within $A$. This third point of view applies more generally  as follows. 

 \begin{lemma}\label{lemdeM} Let $A=\F_2[x]/\<f(x)\>$ for some relation $f(x)=0$. Then we can identify $\bar A\cong\F_2[y]/\<f(1+y)\>=A$ as a change of variables $y=1+x$.
 \end{lemma}\proof  We check this for $f(x)=f_3 x^3+f_2 x^2+ f_1 x+ f_0$. Then $g(y):=f(1+y)=f_3(1+y+y^2+y^3)+f_2(1+y^2)+f_1(1+y)+f_0=f_3y^3 + (f_3+f_2)y^2+ (f_3+f_1)y+ (f_2+f_1+f_0)$.  Hence starting in $\bar A$ and using that $\bar+\mu a=+\mu(1+a)$ in $A$,
 \begin{align*} g_{\bar A}(x)=\bar 0\bar+f_3&x\bar\cdot x^2\bar+ (f_3+f_2) x^2\bar+ (f_3+f_1)x\bar+ (f_3+f_2+f_1+f_0)\bar 1\\
&=
1+f_3x\bar\cdot x^2 + (f_3+f_2) x^2+ (f_3+f_1)x+ (f_3+f_2+f_1+f_0)0 \\
&\quad + f_3+f_3+f_2+f_3+f_1+f_3+f_2+f_1+f_0\\
&=1+ f_3(x^3+x+x^2)+ (f_3+f_2) x^2+(f_3+f_1)x+f_0 \\ &
=1+f_3 x^3+ f_2 x^2+ f_1 x+ f_0=1+f(x)\end{align*}
as an element of $A$. Hence $g_{\bar A}(x)=\bar 0$ in $\bar A$ is equivalent to $f(x)=0$ in $A$,  which in turn is equivalent to a new variable $y=1+x$ with $g(y)=0$ in $A$. \endproof

We also note in passing that for any algebra over $\F_2$, we can define a `generalised derivation'
\[ \del a= a\bar a=a+a^2;\quad \del(ab)=\del a+ \del b+ (\del a)(\del b)=( \del a)\bar\cdot(\del b),\quad\del(a+b)=\del a+\del b.\]
This is just the `infinitesimal part' of the canonical Frobenius automorphism in the sense that the latter is  $F=\id+\del$. Here $\del=0$ for a Boolean algebra, but for a more general algebra we think of it in the spirit of \cite{Law} as a kind of `boundary' of $a$ (the intersection of a subset and its complement, which in the Venn diagram would be the boundary). This not the same as our exterior derivative but is a little similar,  without needing a graph. 

Now let $(\Omega^1,\extd)$ be a differential calculus on $A$ and $\theta\in\Omega^1$. We define $\bar \Omega^1$ as the same set as $\Omega^1$.  

\begin{proposition} Let $(\Omega^1,\extd)$ be a differential calculus on $A$ and  $\theta\in \Omega^1$.  Then $\bar\Omega^1$ defined as the same set as $\Omega^1$ but with a new addition, bimodule structure and differential
\[ \omega\bar+ \eta= \theta+\omega+\eta,\quad a\bar\cdot \omega=a\theta+(a+1)\omega,\quad \omega\bar\cdot a=\theta a+\omega(a+1),\quad  \bar\extd a= \theta+ \extd a\]
is a differential calculus on $\bar A$. Moreover, 

\begin{enumerate}\item[(i)] $\bar{\ }:\Omega^1\to \bar\Omega^1$ defined by $\bar\omega=\theta+\omega$ makes $\bar{\ }:A\to \bar A$ a diffeomorphism. 

\item[(ii)] $\theta$ is the zero element of $\bar\Omega^1$.

\item[(iii)] $\theta$ makes $\Omega^1$ inner if and only if the zero element of $\Omega^1$ makes $\bar\Omega^1$ inner. 
\end{enumerate}
\end{proposition}
\proof Clearly $\omega\bar+(\eta\bar+\zeta)=\omega+\eta+\zeta$ is associative. Moreover
\begin{align*}a\bar\cdot(\omega\bar+\eta)&=a\bar\cdot(\theta+\omega+\eta)=a\theta+ (1+a)\theta+ (1+a)(\omega+\eta) \\
&=\theta+ (\theta+(1+a)\omega)+(\theta+(1+a)\eta)=a\bar\cdot \omega+ a\bar\cdot\eta\\
a\bar\cdot(b\bar\cdot \omega)&=a\theta+ (a+1)(b\bar\cdot\omega)=a\theta+(a+1)b\theta+ (a+1)(b+1)\omega\\ 
&=(a\bar\cdot b)\theta+ (1+(a\bar\cdot b))\omega=(a\bar\cdot b)\bar\cdot\omega\\
a\bar\cdot(\omega\bar\cdot b)&=a\theta+ (a+1)\omega\bar\cdot b=a\theta+(a+1)\theta b+ (a+1)\omega(b+1)\\
&=a\theta+\theta b+ a\theta b+\omega+ a\omega+\omega b+a\omega b=\cdots=(a\bar\cdot \omega)\bar\cdot b, 
\end{align*} 
where we make the same steps in reverse. So $\bar\Omega^1$ is a bimodule. We also have
\begin{align*}(\bar\extd a)\bar\cdot b\bar+ a\bar\cdot\bar\extd b&=\theta+(\theta b+ (\theta+\extd a)(b+1))+(a\theta+ (a+1)(\theta+\extd b))\\
&=\theta+\extd a+\extd b+ (\extd a)b+a\extd b= \bar\extd(a+b+ab)=\bar\extd(a\bar+b). 
\end{align*}
One can check that the surjectivity condition for a differential calculus holds, as it does for $\Omega^1$. 

For the additional facts: (i) Clearly $\bar a\, \bar\cdot\, \bar\omega=(1+a)\bar\cdot(\theta+\omega)=(1+a)\theta+ a(\theta+\omega)=\theta+a\omega=\overline{a\omega}$ and similarly on the other side, so $\bar{\ }:\Omega^1\to\bar\Omega^1$ is a bimodule map in the required sense. The diagram with $\extd,\bar\extd$ also clearly commutes.  (ii)  $\theta\bar+\omega=\theta+\theta+\omega=\omega$ so $\theta$ is the zero element of $\bar\Omega^1$ (iii) $a\, \bar\cdot\,  0=a\theta+(a+1)0=a\theta$ and $0\, \bar\cdot\, a=\theta a$ similarly. Thus $a\, \bar\cdot\, 0\bar+ 0\, \bar\cdot\, a=\theta+[\theta,a]=\bar\extd a$. \endproof

In the example of $\bar P(X)$, we can take $\theta={\rm Arr}$ to have $\bar\Omega^1=\bar P({\rm Arr})$ as a canonical choice, i.e., we recover the procedure in Section~\ref{secDM}. The other canonical choice is $\theta=0$ in which case $\bar\Omega^1$ has an unchanged addition law but a modified product $a\, \bar\cdot\, \omega=\omega\oplus a\cap\omega=\bar a\cap\omega$, the set of arrows in $\omega$ with tip not in $a$. This is not as natural as our previous choice, so we will stick with that. In that case,  we also have $\theta\tens_A\theta={\rm Arr}^{(2)}$ for $P(X)$ which motivates us to similarly define `complementation' on tensor products. With $\omega,\eta\in \Omega^1$ viewed in $\bar\Omega^1$, we define
\[ \omega\tens_{\bar A}\eta:=\omega\tens_A\eta+\theta\tens_A\eta+\omega\tens_A\theta\in \bar\Omega^1\tens_{\bar A}\bar\Omega^1\]
and one can check that $\omega\, \bar\cdot\, a\tens_{\bar A}\eta=\omega\tens_{\bar A}a\, \bar\cdot\,  \eta$. Here $\bar\Omega^1\tens_{\bar A}\bar\Omega^1$ is the same vector space as  $\Omega^1\tens_A\Omega^1$ and, similarly to our treatment of $\bar\Omega^1$,  is a bimodule with 
\begin{equation}\label{barOmega2} \omega\bar+\eta= \theta\tens_A\theta+\omega+\eta,\quad a\, \bar\cdot\, \omega=a\theta\tens_A\theta+(a+1)\omega,\quad \omega\, \bar\cdot\,  a=\theta\tens_A\theta a+\omega(a+1),\end{equation}
where now $\omega,\eta\in \Omega^1\tens_A\Omega^1$. One 
 can check that $\tens_{\bar A}$ is bilinear with respect to $\bar+$. By construction,  we can now define 
\[ \bar{\ }: \Omega^1\tens_{A}\Omega^1\to \bar\Omega^1\tens_{\bar A}\bar\Omega^1,\quad \overline{\omega\tens_A\eta}=\theta\tens_A\theta+ \omega\tens_A\eta\]
and check that it connects the bimodule structures on the two sides compatibly with $\bar{\ }$ on $A$, e.g. on one side this is 
\[ \bar a\, \bar\cdot\, \bar\omega=(a+1)\bar\cdot(\theta\tens_A\theta+\omega)=(a+1)\theta\tens_A\theta+a(\theta\tens_A\theta+\omega)=\theta\tens_A\theta+a\omega=\overline{a\omega}.\]
We now ask when this descends to the wedge product. 

\begin{lemma}  Let $\Omega^1$ extend to an exterior algebra over $A$ at least to $\Omega^2$ and $\extd\theta=0$. 

{\rm (i)}  $\bar\Omega^2$ defined as the same vector space as $\Omega^2$ with bimodule structure as in (\ref{barOmega2}) but now using $\theta^2$ and with $\bar\extd\omega=\theta^2+\extd \omega$ for $\omega\in \Omega^1$ forms the degree 2 part of an exterior algebra over $\bar A$. 

{\rm (ii)}  $\bar{\ }:\Omega^2\to \bar\Omega^2$ defined by $\bar\omega=\theta^2+\omega$ is a map of DGAs to degree 2. 
\end{lemma}
\proof The structure of $\bar\Omega^2$ follows the same structure and proofs as  $\bar\Omega^1\tens_{\bar A}\bar\Omega^1$, namely
\[ \omega\bar+\eta=\theta^2+\omega+\eta,\quad a\bar\cdot\omega=a\theta^2+(a+1)\omega,\quad \omega\bar\cdot a=\theta^2a +\omega(a+1)\]
for $\omega,\eta\in \Omega^2$ and we also have $\omega\bar\cdot\eta=\omega\eta+\theta\eta+\omega\theta$ for $\omega,\eta\in \Omega^1$. We check the Leibniz rule
\begin{align*} \bar\extd(a\bar\cdot\omega)&=\theta^2+\extd(a\theta+(a+1)\omega)=\theta^2+(\extd a)\theta+a\extd\theta+(\extd a)\omega+(a+1)\extd\omega\\
(\bar\extd a)\bar\cdot\omega\bar+a\bar\cdot\bar\extd\omega&=\theta^2+(\theta+\extd a)\omega+\theta\omega+(\theta+\extd a)\theta+a\theta^2+(a+1)(\theta^2+\extd\omega)
\end{align*}
which agree provided $\extd\theta=0$. This is also needed for $\bar\extd\bar\extd a=\theta^2+\extd(\theta+\extd a)=\theta^2+\extd\theta+\extd^2a=\theta^2$ which is the zero element of $\bar\Omega^2$. That $\bar{\ }$ extends our previous map $\Omega^1\to \bar\Omega^1$ compatibly with $\extd$ is also immediate provided $\extd\theta=0$. We also have $\bar\omega\, \bar\cdot\, \bar\eta=(\theta+\omega)(\theta+\eta)+\theta(\theta+\eta)+(\theta+\omega)\theta=\theta^2+\omega\eta=\overline{\omega\eta}$ by construction. 
 \endproof
 
Here $\extd\theta=\{\theta,\theta\}=2\theta^2=0$ is automatic over $\F_2$ if the calculus is inner by $\theta$. This is the case for $P(X)$ with  $\Omega^2_{min}$,  where indeed $\theta\tens_{P(X)}\theta={\rm Arr}^{(2)}$ is the union of all the elements of $\CN_{min}$. We also note that the lemma works similarly for forms of all degrees; we have focussed on the degree 2 case as this is all that is needed for Riemannian geometry.

\subsection{Example of group algebra $A=\F_2\Z_3$.} To illustrate the above,  we focus on the Hopf algebra dual model to the Boolean algebra $\F_2(\Z_3)$ studied in Section~\ref{sectri}. Here, $A=\F_2\Z_3$ is the D model in \cite{MaPac} except that we change $z,x$ there to $x,y$. Then $A$ has basis $1,x,x^2$, the relation $x^3=1$ and the universal calculus $\Omega^1=\Omega^1_{uni}$. 

The new result in this section is to rework the computer results for this algebra in \cite{MaPac} in terms of a left invariant basis $e^+=x^2\extd x$ and $e^-=x\extd x^2$ in a similar spirit to our treatment of $\F_2(\Z_3)$. After a short calculation, the exterior algebra in \cite{MaPac} amounts to the $e^\pm$ as generators and the relations, volume form and inner generator 
\[ e^+x=x(e^++e^-),\quad e^+x^2=x^2e^-,\quad e^-x=x e^+,\quad e^-x^2=x^2(e^++e^-)\]
\[(e^+)^2=(e^-)^2=e^+e^-+e^-e^+=0,\quad{\rm Vol}=e^+e^-,\quad \theta=e^++e^-.\]
In these terms, there are three quantum metrics 
\[ g_i=x^i(e^+\tens e^-+e^-\tens e^+),\quad i=0,1,2\]
(denoted $g_{D.3}, g_{D.1}, g_{D.2}$  in \cite{MaPac}) and they each have a flat QLC
\[ g_0:\quad \nabla e^\pm=0,\quad g_1:\quad \nabla e^+=e^+\tens e^++ g_0,\quad \nabla e^-=e^+\tens e^++ e^-\tens e^-\]
\[g_2:\quad\nabla e^+=e^+\tens e^++e^-\tens e^-,\quad \nabla e^-=e^-\tens e^-+ g_0.\]
Next, each metric has three nonflat equal-curvature connections, with joint curvatures respectively
\[ g_0:\quad R_\nabla e^\pm={\rm Vol}\tens e^\pm,\quad g_1:\quad R_{\nabla}e^+=x{\rm Vol}\tens(e^++e^-),\quad R_{\nabla} e^-=x{\rm Vol}\tens  e^+\]
\[ g_2:\quad R_{\nabla}e^+=x^2{\rm Vol}\tens e^-,\quad R_{\nabla}e^-=x^2{\rm Vol}\tens(e^++e^-)\]
(noting that ${\rm Vol}$ in \cite{MaPac} is $x^2{\rm Vol}$ now). Similarly to the $\F_2(\Z_3)$ model,  \cite{MaPac} tells us that there are two natural lifts that result in an Einstein tensor ${\rm Eins}={\rm Ricci}+g$ that is conserved in the sense $\nabla\cdot {\rm Eins}=0$. When converted to our left-invariant basis, these are
\[ i_\pm({\rm Vol})=e^+\tens e^++ e^-\tens  e^-+e^\pm\tens e^\mp,\]
\[ g_0:\quad  {\rm Eins}_\pm=i_\pm({\rm Vol}),\quad g_1:\quad   {\rm Eins}_\pm= \begin{cases}x\, i_-({\rm Vol})\\ 0,&\end{cases}\  g_2:\quad  {\rm Eins}_\pm= \begin{cases}0\\ x^2 i_+({\rm Vol}).\end{cases}\]
For completeness we also give the three underlying equally curved QLCs from \cite{MaPac} for each metric but converted in terms of our left-invariant forms,
\begin{align*}g_0:\quad &{\rm (i)}\quad\nabla e^+=e^-\tens e^-+g_0,\quad \nabla e^-=e^+\tens e^++e^-\tens e^-\\
&{\rm (ii)}\quad\nabla e^+=e^+\tens e^++e^-\tens e^-,\quad \nabla e^-= e^+\tens e^++g_0\\
& {\rm (iii)}\quad\nabla e^+=e^+\tens e^++g_0,\quad \nabla e^-= e^-\tens e^-+g_0 \\
g_1:\quad &{\rm (i)}\quad\nabla e^+=(1+x^2)e^+\tens e^++x^2 e^-\tens e^-+g_0,\quad \nabla e^-=(1+x^2)e^+\tens e^++e^-\tens e^-+x^2 g_0\\
&{\rm (ii)}\quad\nabla e^+=(1+x^2)(e^+\tens e^++g_0),\quad \nabla e^-= e^+\tens e^++(1+x^2)e^-\tens e^-+x^2g_0\\
& {\rm (iii)}\quad\nabla e^+=e^+\tens e^++x^2 e^-\tens e^-+(1+x^2)g_0    ,\quad \nabla e^-=  (1+x^2)(e^+\tens e^++ e^-\tens e^-)  \\
g_2:\quad &{\rm (i)}\quad\nabla e^+= e^+\tens e^++(1+x)e^-\tens e^-+xg_0 ,\quad \nabla e^-= xe^+\tens e^++(1+x)e^-\tens e^-+g_0\\
&{\rm (ii)}\quad\nabla e^+=  (1+x)e^+\tens e^++e^-\tens e^-+xg_0 ,\quad \nabla e^-=(1+x)(e^-\tens e^-+g_0)   \\
& {\rm (iii)}\quad\nabla e^+= (1+x)( e^+\tens e^++e^-\tens e^-) ,\quad \nabla e^-= xe^+\tens e^++e^-\tens e^-+(1+x)g_0.  \end{align*}
Swapping $e^\pm$ and $x,x^{-1}$ interchanges the $g_1$ and $g_2$ solutions while the $g_0$ solutions transform among themselves with (iii) invariant. In summary, the quantum geometry for the base metric $g_0$ is very similar to that for the Boolean algebra on 3 points in Section~\ref{sec3pt} except that now we have one flat and 3 curved QLCs rather than the other way around, and we also have the possibility of conformally scaled metrics $g_1,g_2$ with slightly different curvatures.

Next, the de Morgan dual algebra by Lemma~\ref{lemdeM} is isomorphic but with dual generator $y= 1+x$ with a new product $x\bar\cdot x=x^2$, $x\bar\cdot x^2=1+x+x^2=x\bar+x^2$. So $x\bar\cdot x^2\bar+ x^2\bar+ x=\bar 0$ in $\bar A$. We also have
$(1+x)^3+(1+x)^2+(1+x)=x+x^2+1+x^2+1+x=0$ so that de Morgan duality is equivalent to a change of variable to $y$ in the same algebra. The associated `derivation' is $\del x=\del x^2=x+x^2$ so that $\del\del x=0$. The de Morgan duality  isomorphism extends to $\Omega^1\cong\bar\Omega^1$, for example
\[ (\bar\extd x)\bar\cdot x\bar+x\bar\cdot \bar\extd x= (\theta+(1+x)\extd x)\bar+(\theta+(\extd x)(1+x))=\theta+\extd(x^2)=\bar\extd(x^2)=\bar\extd(x\bar\cdot x)\]
and this is also necessarily true for $y$ in $\Omega^1$ as $(\extd y)y+y(\extd y)=\extd(y^2)$. 
At degree 2, we have the same vector space for $\bar\Omega^2$  as $\theta^2=0$ and for example
\[ \bar\extd x\bar\cdot\bar\extd (x\bar\cdot x)=(\theta+\extd x)(\theta+\extd x^2)+\theta(\theta+\extd x^2)+(\theta+\extd x)\theta=(\extd x)\extd x^2=x e^+ x^2 e^-=(e^-)^2=0\]
is parallel to  $(\extd y)\extd y^2=0$ in $\Omega$. We can also write
\[ e^-=\overline{e^+}=\theta+x^2\extd x=\theta+(\extd x^2)x=\theta(1+x)+(\theta+\extd x^2)x=(\bar\extd x^2)\bar\cdot(1+x)\]
\[e^+= \overline{e^-}=\theta+x\extd x^2=\theta+(\extd x)x^2=\theta(1+x^2)+(\theta+\extd x^2)x^2=(\bar\extd x)\bar\cdot(1+x^2)\]
as elements of $\bar\Omega^1$, obeying $e^+\bar+e^-=0=\bar\theta$ and 
\[   e^\pm\, \bar\cdot\, e^\pm=\theta e^\pm+e^\pm\theta=0,\quad e^\pm\, \bar\cdot\,  e^\mp=e^\pm e^\mp+(e^++e^-)e^\mp + e^\pm(e^++e^-)=e^\pm e^\mp=e^\mp\, \bar\cdot\, e^\pm.\]
We have $\bar\extd=\extd$ acting on degree 1 and $\overline{\rm Vol}={\rm Vol}$ as $\theta^2=0$, and one can check that the zero element of $\Omega^1$ makes $\bar\Omega$ inner. In short,  $\bar\Omega$ looks different but can also be viewed within $\Omega$ as a change of variables, i.e. de Morgan duality invariance is ultimately part of diffeomorphism invariance. 

We also have, say for the symmetric case (iii) QLC for $g_0$, 
\[ \bar \nabla(e^\pm)=\overline{\nabla{\overline{e^\pm}}}=\overline{\nabla e^\mp}=\overline{e^\mp\tens e^\mp+e^+\tens e^-+e^-\tens e^+}=e^\pm\tens e^\pm\]
which is in the spirit of the curved QLC for $\F_2(\Z_3)$ in Proposition~\ref{proptri}. Equivalently, $\nabla e^\pm$ for the case (iii) QLC is of a similar flavour to the de Morgan dual of this curved QLC in the form stated after Lemma~\ref{lemdemqlctri}. Similarly,  the curvature for $g_0$ has the same $R_\nabla={\rm Vol}\tens$ form as for the curved QLC for the triangle as for its de Morgan dual model in Section~\ref{secdem}. Thus, the Hopf algebra dual model for $\F_2(\Z_3)$ and de Morgan dual model, while very different, also have some striking similarities.  It was  also striking that both the models have  4 QLC's for each metric, just in one case 3 are flat and one is curved and in the other case the situation is reversed. Although $\F_2\Z_3$ has three metrics, these just differ by a scale multiple.  Moreover, in both models there are two natural lifts maps $i$ such that the Einstein tensor is conserved in the sense $\nabla\cdot{\rm Eins}:=((\ ,\ )\tens\id)\nabla({\rm Eins})=0$.

\section{Concluding remarks}\label{seccon}

There are several directions for further work, which we discuss here. The results of Section~\ref{secex} suggest that there is indeed a rich vein of discrete quantum Riemannian geometries on any graph and an immediate task would be to study these systematically for all connected graphs with $|X|=4$ and beyond. The classification of graphs is an unsolved problem and our results suggest a new geometric class of invariants based on the quantum Riemannian geometries that they support. One could also consider geometric invariants more broadly, such as could be obtained from quantum geometric Chern-Simons theory and related TQFTs on a graph. These questions apply over any field and relate to other efforts in noncommutative algebra, such as the recent notion of a Hopf algebroid of differential operators\cite{Gho} reconstructed from the moduli of flat bimodule connections more broadly (not necessarily on $\Omega^1$ as studied here). The digital case over $\F_2$ should be seen as an extreme case where the moduli spaces reduce more directly to combinatorics with calculable results, as we have seen here by algebraic means and previously in \cite{MaPac} by computer means. We are also not obliged to look at the discrete set case and could look both at Boolean algebras more generally\cite{Hai} than we have done and at other types of algebras over $\F_2$. Examples of the latter,  where the quantum Riemannian geometry remains to be explored, include $M_2(\F_2)$ using the classification of $\Omega^1$ in \cite{BegMa}, the commutative and cocommutative Hopf algebras $A_d$ in \cite{BasMa}, and the noncommutative and noncocommutative Hopf algebra ${\sl dsl}_2$ in \cite{MaPac1}, where bicovariance leads to a natural construction for exterior algebras $\Omega$.  TQFT in some form, such as the Kitaev model, applies to any Hopf algebra\cite{Meu} but their $\F_2$-versions also remain to be explored.

Our other main theme was the Boolean idea of de Morgan duality, which we showed extends to the quantum Riemannian geometry in the complete  atomic Boolean case of the power set $P(X)$. That quantum Riemannian geometry produces reasonable answers in the form of $\cup$ and $\cap$ and Venn diagrams now including graphs speaks to the robustness of the formalism and opens up a self-contained set-up that can be explored further and with more attention to the case of infinite $X$, as mentioned at the end of Section~\ref{secboo}. It is, however, fair to say that the actual applications to logic and the significance of curvature there remain to be explored. In everyday life, the `logic' of subsets of a set does not necessarily make reference to a graph structure on $X$, but it can do in the context of some kind of process where one element can turn into another. A general class of interest in computer science would be $(X,\le)$ a preorder (a transitively complete graph extended to include all self-arrows, where $x\to y$ means $x\le y$). These could be used, for example, to encode chemical or manufacturing processes\cite{Spiv}. For another example, if $Y$ is a set without structure then $X=P(Y)$ is a preorder by subset inclusion, or in  propositional logic terms the extended graph arrows are implication $\Rightarrow$. In this context, a connection could allow parallel transport of elements of proofs and curvature could potentially acquire an interpretation. However, the directed graph here would not be bidirectional and one would either need to work with degenerate metrics or look at other bundles than the cotangent bundle. 

We  saw in particular that the de Morgan dual of a differential form makes sense and amounts to complementation in the set of arrows. For example, we saw that the element $\theta={\rm Arr}$ consisting of all arrows maps to the zero differential form or empty set of arrows and that this indeed is the inner element for the de Morgan dual differential calculus. We have not attempted to discuss any physics but one could say in some sense that a differential form of `maximum density' (all arrows switched on) maps over to one of `zero density' (all arrows switched off) somewhat in the spirit of some kind of gravity - quantum duality\cite{Ma:pri,Ma:qg,Ma:mar}. What we saw was that  such a duality can nevertheless be formulated precisely as a diffeomorphism between $P(X)$ and the dual model on $\bar P(X)$ where $\cup,\cap$ are swapped. This is a map between different algebras but on the same set, so it could be viewed as some kind of `coordinate transform', except that this is not simply a change of generators as the algebra structure is being changed. Nevertheless, if two manifolds are diffeomorphic then Riemannian geometry on one is equivalent to Riemannian geometry on the other, and our interpretation was somewhat analogous to this. 

It is an interesting question as to how de Morgan duality then generalises, and in Section~\ref{secbarA} we gave one answer from the point of view of $\F_2$-algebras more generally, and observed some similarities with Hopf algebra duality. The latter has already been proposed in quantum gravity as `quantum Born reciprocity' and is somewhat different in character from de Morgan duality. Its interaction with quantum Riemannian geometry was previously considered in  \cite{MaTao2} between  $\C(S_3)$ and $\C S_3$, and potentially other finite groups. In the first model, the possible differential structures are labelled by conjugacy classes and the eigenvectors of the resulting Laplacians `waves' provided by matrix elements of irreducible representations, in the dual model the possible differential structures are labelled by representations and eigenvectors of the Laplacians by conjugacy classes. Over $\C$, the two $\Z_3$ models would be isomorphic by Hopf algebra self-duality of Fourier transform, but this is not the case over $\F_2$. 

Another direction for generalisation would be from Boolean algebras to weaker structures of interest in logic, such as Heyting algebras and lattices and it could be interesting to develop quantum Riemannian geometry for these, which could be of interest both in topos theory and possibly in the foundations of quantum mechanics\cite{DorIsh}.  It is also possible for a Heyting algebra to obey the de Morgan duality identities relating $\cap$ and $\cup$ even if $a\cup\bar a\ne 1$, but with double negation no longer the identity if we want to be beyond the Boolean case. Heyting algebras in fact arise in many contexts and a natural example is $C(X,[0,1])$ where $X$ is a discrete set and we replace $\{0,1\}$ for the Boolean case by `probabilities' with values between $[0,1]$. The $\cap, \cup$ (meet and join) are given pointwise as the min, max of the values and the Heyting negation of a function $f$ is the characteristic function of its zero set. Here $(f\cap\bar f)(x)=\min(f(x),\bar f(x))=0$ at all $x$ but $(f\cup\bar f)(x)=\max(f(x),\bar f(x))\ne 1$, having value $1$ where $f(x)=0$ but $f(x)$ otherwise. In this case, we can still view $\bar{\ }$ as an algebra map to the dual structure in the same spirit as our treatment in the Boolean case, and ask about the extension to differentials. The problem here is two fold: first, $\bar{\ }$ is no longer an isomorphism and, second, $\cup$ does not lead to a proper addition law; we do not in fact have an algebra exactly and must therefore also approach the Leibniz rule differently. In the particular case of $C(X,[0,1])$, the first problem can be solved by using a different  `complementation' $\bar f(x)=1-f(x)$ which now does square to the identity and also interchanges $\cup$ and $\cap$ in a de Morgan like manner (so gives an isomorphism to the $\cup,\cap$-reversed algebra). One also has a second product, namely the usual pointwise product of functions, and a candidate `addition' $f+g-fg$ but neither product distributes over it. Another approach that could be considered here is that of  bi-Heyting algebras\cite{ReyZol}. A probabilistic setting would also connect to the idea that metrics in quantum Riemannian geometry do not need to be edge symmetric (there could be a different `length' in the two directions for each edge). In \cite{Ma:mar}, 
we have proposed instead that such asymmetric edge weights be Markov transition probabilities with values in $[0,1]$, giving a quantum-geometric picture of Markov processes on stochastic vectors $f\in C(X,[0,1])$ viewed within the actual algebra $C(X,\R)$. De-Morgan duality here remains to be considered but could ultimately re-emerge in a probabilistic interpretation. This is far from the $\F_2=\{0,1\}$ valued models in the present paper but could be seen as a natural generalisation where we allow intermediate values. These are some ongoing directions for further work.

\end{document}